\documentclass[a4paper,10.5pt]{article}
\usepackage{amsmath,amssymb,amsthm}
\theoremstyle{definition}
 \newtheorem{dfn}{Definition}[section]
 \newtheorem{remark}[dfn]{Remark}

\theoremstyle{plain}
 \newtheorem{thm}[dfn]{Theorem}
 
 \newtheorem{lem}[dfn]{Lemma}
 \newtheorem{cor}[dfn]{Corollary}

\numberwithin{equation}{section}
\newcommand{\bn}{{\bold n}}
\newcommand{\bff}{{\bold f}}

\newcommand{\bD}{{\bold D}}

\newcommand{\bF}{{\bold F}}
\newcommand{\bG}{{\bold G}}
\newcommand{\bH}{{\bold H}}
\newcommand{\bI}{{\bold I}}
\newcommand{\bJ}{{\bold J}}
\newcommand{\bK}{{\bold K}}
\newcommand{\bL}{{\bold L}}

\newcommand{\bT}{{\bold T}}

\newcommand{\bV}{{\bold V}}
\newcommand{\bW}{{\bold W}}
\newcommand{\bX}{{\bold X}}
\newcommand{\DV}{{\rm Div}\,}
\newcommand{\dv}{{\rm div}\,}

\newcommand{\ba}{\bold a}
\newcommand{\bb}{\bold b}
\newcommand{\bx}{\bold x}
\newcommand{\bv}{\bold v}
\newcommand{\bu}{\bold u}
\newcommand{\bw}{\bold w}
\newcommand{\bz}{\bold z}

\newcommand{\bg}{\bold g}
\newcommand{\bh}{\bold h}
\newcommand{\bp}{\bold p}
\newcommand{\bZ}{\bold Z}

\newcommand{\BR}{{\Bbb R}}
\newcommand{\BC}{{\Bbb C}}

\newcommand{\BN}{{\Bbb N}}

\newcommand{\BI}{{\Bbb I}}
\newcommand{\BJ}{{\Bbb J}}
\newcommand{\BK}{{\Bbb K}}

\newcommand{\BT}{{\Bbb T}}

\newcommand{\BZ}{{\Bbb Z}}

\newcommand{\CA}{{\mathcal A}}

\newcommand{\CC}{{\mathcal C}}
\newcommand{\CD}{{\mathcal D}}

\newcommand{\CF}{{\mathcal F}}
\newcommand{\CI}{{\mathcal I}}

\newcommand{\CL}{{\mathcal L}}

\newcommand{\CR}{{\mathcal R}}
\newcommand{\CS}{{\mathcal S}}

\newcommand{\CW}{{\mathcal W}}
\newcommand{\CX}{{\mathcal X}}

\newcommand{\pd}{\partial}

\newcommand{\Hol}{{\rm Hol}}

\begin{document}
\title{On some free boundary problem of the Navier-Stokes equations\\
in the maximal $L_p$-$L_q$ regularity class}
\author{Yoshihiro SHIBATA
\thanks{Partially supported by JST CREST and JSPS
Grant-in-aid for Scientific Research (S) \# 24224004 }
\\
{\small Department of Mathematics and Research Institute of 
Science and Engineering} \\
{\small Waseda University, 
Ohkubo 3-4-1, Shinjuku-ku, Tokyo 169-8555, Japan} \\
{\small e-mail address: yshibata@waseda.jp}
}
\date{}
\maketitle
\begin{abstract}
This paper is concerned with the free boundary problem for
the Navier Stokes equations without surface tension
in the $L_p$ in time and $L_q$ in space setting
with $2 < p < \infty$ and $N < q < \infty$. 
A local in time existence theorem is proved in a 
uniform $W^{2-1/q}$ domain in the $N$-dimensional Euclidean
space $\BR^N$ ($N \geq 2$) under the assumption that
weak Dirichlet-Neumann problem is uniquely solvable. 
Moreover, a global in time existence theorem is 
proved for small initial data 
under the assumption that $\Omega$ is bounded 
additionally.  This was already proved by Solonnikov
\cite{Sol1} by using the continuation argument of 
local in time solutions which are exponentially
stable in the energy level under the assumption that the initial 
data is orthogonal to the rigid motion. 
 We also use the continuation argument and the same orthogonality
for the initial data. But,  our argument about the continuation of 
local in time solutions  is based on 
some decay theorem for the linearized problem,
which is a different point than \cite{Sol1}. 
\end{abstract}
\section{Introduction} \label{sec:1}
The present paper deals with some local and global in time 
unique existence theorems of solutions to the Navier-Stokes equations
describing the motion of  a viscous
incompressible fluid flow with free surface 
without taking surface tension into account.  
Our problem is formulated 
in the following.  Let $\Omega$ be a domain in the 
$N$-dimensional Euclidean space
$\BR^N$ ($N \geq 2$)
occupied by a viscous incompressible fluid.  We assume that 
the boundary of $\Omega$ consists of two parts $S$ and $\Gamma$ with
$S \cap \Gamma = \emptyset$. We may assume that 
$\Gamma$ is an empty set. Let $\Omega_t$ and 
$S_t$ be evolutions of $\Omega$ and $S$ with time 
variable $t > 0$ and we assume that $S_t \cap \Gamma=\emptyset$ for $t\geq 0$. 
The velocity vector field $\bv = 
\bv(x, t) = (v_1(x, t), \ldots, v_N(x, t))$ and the pressure
$\pi=\pi(x,t)$ for $x =(x_1, \ldots, x_N) \in \Omega_t$ 
satisfy the Navier-Stokes equations
\begin{equation}\label{1.1} 
(\pd_t\bv + (\bv\cdot\nabla)\bv) -\DV \bT(\bv, \pi) = 0, 
\quad \dv \bv = 0.
\end{equation}
The initial conditions, the boundary conditions on the free boundary
$S_t$ and the non-slip conditions on the fixed boundary $\Gamma$
have the following forms: 
\begin{equation}\label{1.2}\begin{split}
\bv|_{t=0} &= \bv_0
\quad\text{in $\Omega$}, \\ 
\bT(\bv, \pi)\bn_t|_{S_t} &= 0,
\quad \bv|_{\Gamma} = 0.
\end{split}\end{equation}
Here, 
$\bn_t$ is the unit outward normal to $S_t$.
Moreover, $\bT = \bT(\bv, \pi)$ denotes the stress tensor of the form:
\begin{equation}\label{1.3}
\bT(\bv, \pi) = -\pi\bI + \mu\bD(\bv) 
\end{equation}
where $\mu$ denotes a positive constant describing the  
viscosity coefficient,
$\bD(\bv)$  the deformation tensor whose $(j, k)$ components
are $D_{jk}(\bv) = (\pd_jv_k + \pd_kv_j)$ with $\pd_j
= \pd/\pd x_j$,  and $\bI$  the $N\times N$ identity matrix. Finally, 
for any matrix field $\bK$ with components $K_{ij}$, $i, j=1, \ldots, N$,
the quantity $\DV \bK$ is an $N$-vector with $i$-th component 
$\sum_{j=1}^N\pd_j K_{ij}$, and also for any vector of functions
$\bu=(u_1, \ldots, u_N)$ we set $\dv \bu =\sum_{j=1}^N \pd_ju_j$,
$\bu\cdot\nabla = \sum_{j=1}^Nu_j\pd_j$ and $\pd_t\bu =
(\pd u_1/\pd t, \ldots, \pd u_N/\pd t)$. 

Aside from the dynamical system \eqref{1.1},  we impose 
a further kinematic condition: 
\begin{equation}\label{kinetic:1}
\pd_tF + (\bv\cdot \nabla) F = 0 \quad\text{on $S_t$},
\end{equation}
where $S_t$ is defined by $F = F(x, t) = 0$ locally. 
In other words, $S_t$ is given by 
\begin{equation}\label{1.4}
S_t = \{x  \in \BR^N \mid x = 
\bx(\xi, t) \enskip (\xi \in S) \}, 
\end{equation}
where $\bx = \bx(\xi, t)$ is the solution to the Cauchy problem:
$
\dot\bx = d\bx/dt= \bv(\bx, t) \enskip(t > 0)$ with
$\bx|_{t=0} = \xi$.  
This expresses the fact that the free boundary $S_t$ consists of 
the same particles for all $t > 0$, 
which do not leave it and are not incident from $\Omega_t$.


The free boundary problem for the Navier-Stokes equations has been 
studied by many mathematicians in the following two cases:
\begin{itemize}
\item[\thetag1]~The motion of an isolated liquid mass;
\item[\thetag2]~The motion of a viscous incompressible fluid
 contained in an ocean of infinite content.
\end{itemize}
In case \thetag1 the initial domain $\Omega$ is bounded. A local in 
time unique existence theorem was proved by Solonnikov
\cite{Sol2, Sol4, Sol5, Sol6} in the $L_2$ Sobolev-Slobodetskii
space, by Schweizer \cite{Sch} in the semigroup setting, 
by Moglievski\u\i\, and Solonnikov \cite{Mog-Sol, Sol6} in the 
H\"older spaces with surface tension; and by Solonnikov 
\cite{Sol1}  and Mucha and W.~Zaj\c aczkowski \cite{Much-Zaj2} in the 
$L_p$ Sobolev-Slobodetskii space and by Shibata and Shimizu 
\cite{SS1} in the $L_p$ in time and $L_q$ in space setting without
surface tension. A global in time unique existence theorem for 
small initial velocity was proved by Solonnikov \cite{Sol1}
in the $L_p$ Sobolev-Slobodetskii space without surface tension; 
and by Solonnikov \cite{Sol3}
in the $L_2$ Sobolev-Slobodetskii space and by Padula and 
Solonnikov \cite{Pad-Sol} in the H\"older spaces under the 
additional assumption that the initial domain $\Omega$ is sufficiently 
close to a ball with surface tension.  

In case \thetag2, the initial domain $\Omega$ is a perturbed layer 
like: $\Omega = \{x \in \BR^N \mid -b < X_N < \eta(x'), 
x' = (x_1, \ldots, x_{N-1}) \in \BR^{N-1}\}$.  
A local in time unique existence theorem was proved by 
Beale \cite{Beale1}, Allain \cite{Allain} and Tani \cite{Tani} in the 
$L_2$ Sobolev-Slobodetskii space with surface tension 
and by Abels \cite{Abels} in the $L_p$ Sobolev-Slobodetskii space
without surface tension. A global in time unique existence
theorem for small initial velocity was proved in the 
$L_2$ Sobolev-Slobodetskii space by Beale \cite{Beale2} and 
Tani and Tanaka \cite{Tani-Tanaka} with surface tension, and 
by Sylvester \cite{Slv1} without surface tension.  The decay
rate was studied by Beale and Nishida \cite{Beale-Nishida},
Sylvestre \cite{Slv2} and Hataya \cite{Hataya}. 

The purpose of this paper is to prove a local in time 
unique existence theorem for problem
\eqref{1.1} and \eqref{1.2} under the assumption
that the initial domain $\Omega$ is a uniform $W^{2-1/q}_q$
($N < q < \infty$) domain and weak Dirichlet-Neumann 
problem is uniquely solvable \footnote{These 
assumptions are exactly stated in Definition \ref{dfn:1.1}
and Definition \ref{dfn:1.2} in the following.}, 
which includes the cases
\thetag{1} and \thetag{2} without surface tension.
And also, we prove  a global in time unique 
existence theorem for problem
\eqref{1.1} and \eqref{1.2}
for a small initial data 
in the $L_p$ in time and $L_q$ in space setting
assuming that $\Omega$ is bounded in addition.
This was mentioned in Shibata and Shimizu \cite{SS1},
but there was a serious gap in the proof, so that
we reprove it in a  different approach than \cite{SS1} 
in this paper. 

To prove a local in time unique existence theorem,
the key step is to prove the maximal regularity theorem 
for the linearized equations given in the following:
\begin{gather}
\pd_t\bu - \DV\bT(\bu, \theta) = \bff, \quad
\dv \bu = g = \dv\bg \quad  \text{in $\Omega\times(0, T)$}, 
\nonumber\\
\bT(\bu, \theta)\tilde\bn|_{S}  = \bh|_{S}, 
\quad \bu|_{\Gamma} = 0, 
\quad
\bu|_{t=0}  = \bu_0 \quad\text{in $\Omega$}
\label{lp}
\end{gather}
with $0 < T \leq \infty$. 
Here, $\tilde\bn$ denotes the extension of $\bn$ to the whole space
$\BR^N$.  In fact, as was seen in \cite[\thetag{5.12}]{EBS} (cf. also
\cite[Appendix]{ES}), we can define $\tilde\bn$ on $\BR^N$ such that
$\tilde\bn|_{S} = \bn$ and   
\begin{equation}\label{ext-normal} 
\|f\tilde\bn\|_{W^1_q(\Omega)} \leq C\|f\|_{W^1_q(\Omega)}
\end{equation}
for any $f \in W^1_q(\Omega)$ with some constant $C$ depending on
$\Omega$ if $\Omega$ is a uniform $W^{2-1/r}_r$ domain with
$N < r < \infty$.

To prove the maximal regularity theorem, problem \eqref{lp} is reduced locally 
to the model problems 
in a neighbourhood of either an interior point or a
boundary point by using the localization technique and 
the partition of unity associated with the domain 
$\Omega$. The boundary neighbourhood problem \eqref{lp} is 
transformed to a problem in the half-space $x_N > 0$.  By applying the 
Fourier transform with respect to time and tangential directions,
problem \eqref{lp} becomes a system of ordinary differential 
equations.  Solonnikov \cite{Sol0} calculates explicityly the 
inverse Fourier transform of solutions of such ordinary 
differential equations and expresses them 
in the form of potentials in the half-space.  Then, he 
estimates them in suitable norms. 
Mucha and Zaj\c aczkowski \cite{Much-Zaj1} directly 
estimate them 
using the multiplier theorem of
Marinkiewicz and Mikhlin type \cite{Mihlin}. 

On the other hand, Shibata \cite{S1} proved the maximal regularity 
theorem\footnote{ The maximal regularity theorems are 
given in Theorem \ref{thm:linear} and Theorem \ref{thm:linear*}
in Sect. 2 in the following.} 
by using the $\CR$- bounded solution operators to the 
corresponding resolvent problem of the form:
\begin{equation}\label{gs}\begin{split}
\lambda \bv - \DV\bT(\bv, \kappa) = \bff, \quad
\dv\bv &= g = \dv\bg \quad\text{in $\Omega$}, \\
\bT(\bv, \kappa)\tilde\bn|_{S} = \bh|_{S}, \quad
\bv|_{\Gamma} & = 0.
\end{split}\end{equation}
In fact, according to the theorem in \cite{S1}, 
for any $\epsilon \in (0, \pi/2)$ 
there exist a constant $\lambda_0 \geq 1$ and an operator family
$\CR(\lambda) \in \Hol(\Sigma_{\epsilon, \lambda_0}, 
\CL(\CX_q(\Omega), W^2_q(\Omega)^N))$ such that for any 
$\bff \in L_q(\Omega)^N$, $g \in W^1_q(\Omega)$, 
$\bg \in L_q(\Omega)^N$ and $\bh \in W^1_q(\Omega)^N$, 
problem \eqref{gs} admits a unique solution 
$\bv = \CR(\lambda)(\bff, \lambda^{1/2}g, \nabla g, \lambda\bg, \lambda^{1/2}
\bh, \nabla\bh)$ with some pressure term $\kappa$, 
and $(\lambda, \lambda^{1/2}\nabla, \nabla^2)\CR(\lambda)$ 
is $\CR$ bounded for $\lambda \in \Sigma_{\epsilon, \lambda_0}$ with 
value in $\CL(\CX_q(\Omega), L_q(\Omega)^{\tilde N})$. Here,  $\tilde N
= N+N^2+N^3$, $\Sigma_{\epsilon, \lambda_0}
= \{\lambda \in \BC \mid |\lambda| \geq \lambda_0, 
|\arg\lambda| \leq \pi-\epsilon\}$,
 $\CX_q(\Omega) = \{F = (F_1, \ldots, F_6) \mid
F_1, F_3, F_4, F_5 \in L_q(\Omega)^N, F_2 \in L_q(\Omega), 
F_6 \in L_q(\Omega)^{N^2}\}$, and $F_1$, $F_2$, $F_3$, $F_4$, $F_5$ and 
$F_6$ are independent variables corresponding to $\bff$, $\lambda^{1/2}g$, 
$\nabla g$, $\lambda\bg$, $\lambda^{1/2}\bh$ and $\nabla\bh$, respectively. 
Moreover, $\Hol(\Sigma_{\epsilon, \lambda_0}, \CL(X, Y))$ denotes the set of 
all $\CL(X, Y)$ valued holomorphic functions defined on 
$\Sigma_{\epsilon, \lambda_0}$ and $\CL(X, Y)$ the set of all bounded 
linear operators from a Banach space $X$ into another Banach space $Y$. 
Since the solution $\bu$ for \eqref{lp}
is given by the Laplace inverse transform
of $\CR(\lambda)(\bff, \lambda^{1/2}g, \nabla g, \lambda\bg, \lambda^{1/2}
\bh, \nabla\bh)$, the maximal regularity is obtained 
with help of Weis' operator valued Fourier multiplier
theorem \cite{Weis}.  

Finally, we introduce some symbols used throughout the paper. 
For any domain $D$ and $1\leq q\leq \infty$, $L_q(D)$ and $W^m_q(D)$  
denote the usual Lebesgue space 
and Sobolev space, while $\|\cdot\|_{L_q(D)}$ and $\|\cdot\|_{W^m_q(D)}$
denote their norms, respectively.  We set $W^0_q(D) = L_q(D)$.
$C^\infty_0(D)$ denotes the set of all 
$C^\infty(\BR^N)$ functions whose supports are compact
and contained in $D$. We set $(f, g)_D = \int_D f(x)g(x)\,dx$.  
For any Banach space $X$ and 
$1 \leq p \leq \infty$, $L_p((a, b), X)$ and 
$W^m_p((a, b), X)$ denote the usual  Lebesgue space and Sobolev space of
$X$-valued functions defined on an interval $(a, b)$, 
while $\|\cdot\|_{L_p((a, b), X)}$ and 
$\|\cdot\|_{W^m_p((a, b), X)}$ denote their norms, respectively. 
For $0 < \theta < 1$,  
$B^{2\theta}_{q,p}(D)$ denotes the real interpolation 
space defined by $B^{2\theta}_{q,p}(D)
= (L_q(D), W^2_q(D))_{\theta, p}$ with 
real interpolation functor $(\cdot, \cdot)_{\theta, p}$, while
$\|\cdot\|_{B^{2\theta}_{q,p}(D)}$ denotes its norm.  We set 
$W^{2\theta}_q = B^{2\theta}_{q,q}$. 
The $d$-product space of $X$ is defined by  
$X^d = \{f = (f, \ldots, f_d) \mid 
f_i \in X (i=1, \ldots, d)\}$, while its norm is denoted by 
$\|\cdot\|_X$ instead of $\|\cdot\|_{X^d}$ for the sake of 
simplicity. 
$\BN$, $\BR$ and $\BC$ denote the sets of all natural numbers,
real numbers and complex
numbers, respectively. We set $\BN_0 = \BN \cup \{0\}$. 
For any multi-index $\kappa = (\kappa_1, \ldots, \kappa_N) 
\in \BN_0^N$, we write $|\kappa| = \kappa_1 + \cdots + \kappa_N$ and
$\pd_x^\kappa = \pd_1^{\kappa_1}\cdots\pd_N^{\kappa_N}$ with 
$x = (x_1, \ldots, x_N)$ and $\pd_j = \pd/\pd x_j$. 
For any scalor function $f$ and $N$-vector of functions 
$\bg$, we set 
\begin{alignat*}2
\nabla f &= (\pd_1f, \ldots, \pd_Nf), &\enskip
\nabla \bg &= (\pd_ig_j \mid i, j=1, \ldots, N), \\
\nabla^2f & = (\pd^\alpha f \mid |\alpha = 2 ),
&\enskip 
\nabla^2\bg & = (\pd^\alpha g_i \mid |\alpha|=2, i=1, \ldots, N).
\end{alignat*}
For $\ba = (a_1, \ldots, a_N)$ and $\bb=(b_1, \ldots, b_N)
\in \BR^N$, we set $\ba\cdot\bb = <\ba, \bb> = 
\sum_{j=1}^Na_jb_j$. For scalor functions $f$, $g$
and $N$-vectors of functions $\bff$, $\bg$, we set 
$(f, g)_D = \int_D f(x)g(x)\,dx$ and $(\bff, \bg)_D
= \int_D\bff(x)\cdot\bg(x)\,dx$. 
The letter $C$ denotes generic constants and 
the constant $C_{a, b, \cdots}$ depends on $a$, $b$, $\cdots$.
The values of constants $C$ and $C_{a, b, \cdots}$ may change from line to 
line. 

\section{Main Results} 
In this section, we state our main results. 
Since $\Omega_t$ should be decided, we transfer $\Omega_t$ to $\Omega$
by  the Lagrange transformation  as follows:  If the 
velocity field $\bu(\xi, t)$ is known as a 
function of the Lagrange coordinates $\xi \in \Omega$, then the Euler 
coordinates $x \in \Omega_t$ is written in the form:
$$
x = \xi + \int^t_0 \bu(\xi, s)\,ds \equiv \bX_\bu(\xi, t),
$$
where $\bu(\xi, t) = (u_1(\xi, t), \ldots, u_N(\xi, t))
= \bv(\bX_\bu(\xi, t), t)$.  Let $A$ be the Jacobi matrix of the 
transformation $x = \bX_\bu(\xi, t)$ with elements $a_{ij} = \delta_{ij}
+ \int^t_0(\pd u_i/\pd \xi_j)(\xi, s)\,ds$.  Since 
$\det A = 1$ as follows from $\dv \bv=0$ in $\Omega_t$, denoting
the cofactor matrix of  $A$ by $\CA$, we have $\nabla_x = \CA\nabla_\xi$
with $\nabla_x ={}^T(\pd/\pd x_1, \ldots, \pd/\pd x_N)$ and 
$\nabla_\xi = {}^T(\pd/\pd\xi_1, \ldots, \pd/\pd\xi_N)$ 
\footnote{${}^TM$ denotes the transposed $M$.}.
We can represent $\CA$ by 
$\CA = \bI + \bV_0(\int^t_0\nabla\bu(\xi, s)\,ds)$ 
with some matrix $\bV_0(\bK)$  of polynomials with respect to
 $\bK = (k_{ij})$ satisfying the condition: $\bV_0(0) = 0$, where 
$k_{ij}$ is a 
corresponding variable to $\int^t_0 (\pd u_i/\pd \xi_j)(\xi, s)\,ds$. 
Let $\bn$ be the unit outward normal to $S$, and then 
by \eqref{kinetic:1} we have 
\begin{equation}\label{1.7}
\bn_t= \frac{\CA\bn}{|\CA\bn|}.
\end{equation}
We also see that 
\begin{equation}\label{1.8}
\dv_x \bw = \dv_\xi({}^T\CA\hat \bw) = {\rm tr}\,(\CA\nabla_\xi\hat \bw)
\end{equation}
with $\hat \bw(\xi, t) = \bw(\bX_\bu(\xi, t), t)$, where ${\rm tr}\,M$
denotes the trace of any matrix $M$. Moreover, 
 what $A=(A^{-1})^{-1}
= \CA^{-1}$ yields that 
\begin{equation}\label{1.9} \CA^{-1} =
\bI + \bV_{1}(\int^t_0\nabla \bu(\xi, s)\,ds)
\end{equation}
with some matrix $\bV_1(\bK)$ of polynomials with respect to
 $\bK = (k_{ij})$ satisfying the condition:$\bV_1(0) = 0$. 
Using \eqref{1.7}, \eqref{1.8} and \eqref{1.9}, and 
setting $\theta(\xi, t) = \pi(\bX_\bu(\xi, t))$, we have the 
following Lagrangian description of problem \eqref{1.1}-\eqref{1.2}:
\begin{gather}
\pd_t\bu - \DV\bT(\bu, \theta) 
= \bF(\bu), \enskip 
\dv \bu = G(\bu)  = \dv \bG(\bu) \quad\text{in $\Omega\times(0, T)$}, 
\nonumber\\
\bT(\bu, \theta)\tilde\bn|_S  = \bH(\bu)\tilde\bn|_S, 
\enskip
\bu|_\Gamma=0, \enskip 
\bu|_{t=0}  = \bv_0 \enskip\text{in $\Omega$}.
\label{p1}
\end{gather}
Here, $\bF(\bu)$, $g(\bu)$, $\bg(\bu)$ and $\bH(\bu)$ are nonlinear
functions of the forms:
\begin{align}
\bF(\bu) & = -\bV_{1}(\int^t_0\nabla\bu\,ds)\pd_t\bu + 
\bV_2(\int^t_0\nabla\bu\,ds)\nabla^2\bu + \bV_3(\int^t_0\nabla\bu\,ds)
\int^t_0\nabla^2\bu\,ds\cdot \nabla\bu, \nonumber\\
G(\bu) & = \bV_4(\int^t_0\nabla\bu\,ds)\nabla\bu, 
\quad \bG(\bu) = \bV_5(\int^t_0\nabla\bu\,ds)\bu,
\quad
\bH(\bu) = \bV_6(\int^t_0\nabla\bu\,ds)\nabla\bu,
\label{nonlinear:1}
\end{align}
with some matrices $\bV_i(\bK)$ $(i = 1, \ldots, 6)$ 
of polynomials with respect to $\bK$ satisfying the conditions: 
\begin{equation}\label{null:1}
\bV_1(0) = 0, 
\bV_2(0)= 0, \enskip \bV_4(0)= 0,\enskip \bV_5(0)= 0, \enskip \bV_6(0)= 0.
\end{equation}


We introduce the definition of uniform $W^{2-1/r}_r$ domain.
\begin{dfn}\label{dfn:1.1}
Let $1 < r < \infty$ and let $\Omega$ be a domain in 
$\BR^N$ with boundary $\pd\Omega$.  We say that  $\Omega$ is a 
uniform $W^{2-1/r}_r$ domain,  
if there exist positive constants 
$\alpha$, $\beta$ and $K$ such that for any $x_0 = (x_{01}, \ldots, 
x_{0N}) \in \pd\Omega$ 
there exist a coordinate number $j$ and
 a $W^{2-1/r}_r$ function $h(x')$ $(x' = (x_1, \ldots, 
\hat x_j, \ldots, x_N))$  
defined on $B'_\alpha(x_0')$ with $x_0' = (x_{01},
\ldots \hat x_{0j}, \ldots, x_{0N})$ and 
$\|h\|_{W^{2-1/r}_r(B'_\alpha(x_0'))}
\leq K$ such that 
\begin{equation}\label{eq:gd}\begin{split}
\Omega\cap B_\beta(x_0) &= \{x \in \BR^N \mid 
x_j > h(x')\enskip (x' \in B'_\alpha(x'_0))\} \cap 
B_\beta(x_0), \\
\pd\Omega\cap B_\beta(x_0) &= \{x \in \BR^N \mid 
x_j = h(x')\enskip (x' \in B'_\alpha(x'_0))\} \cap 
B_\beta(x_0).
\end{split}\end{equation}
Here, $(x_1, \ldots, \hat x_j, \ldots, x_N)
= (x_1, \ldots, x_{j-1}, x_{j+1}, \ldots, x_N)$, 
$B'_\alpha(x_0') = \{x' \in \BR^{N-1} \mid 
|x' - x_0'| < \alpha\}$ and $B_\beta(x_0) = \{x \in \BR^N \mid
|x - x_0| < \beta\}$.
\end{dfn}
To prove our local in time unique existence theorem for \eqref{p1} 
in  a uniform $W^{2-1/q}_q$ domain,
we need the unique solvability of
weak Dirichlet-Neumann problem to treat the divergence condition.
But, in general it is not known except for the $L_2$ framework,
so that we have to assume it in this paper. For this purpose, 
we introduce spaces  
$W^1_{q, 0}(\Omega)$ and $\hat W^1_{q, 0}(\Omega)$
 defined by 
$\hat W^1_{q, 0}(\Omega) = \{ \theta \in L_{q, {\rm loc}}
(\Omega) \mid \nabla\theta \in L_q(\Omega)^N, \theta|_{S} = 0\}$ and  
$W^1_{q, 0}(\Omega) = \{ \theta \in  W^1_q(\Omega) \mid 
\theta|_S = 0\}$. 
\begin{dfn}\label{dfn:1.2} ~ Let $1 < q < \infty$ and 
let $\CW^1_q(\Omega)$ be a closed subspace of $\hat W^1_{q, 0}
(\Omega)$ that contains $W^1_{q, 0}(\Omega)$.  
Then,  weak Dirichlet-Neumann problem is called 
uniquely solvable for $\CW^1_q(\Omega)$,  
 if the following assertion holds:  For any
$f \in L_q(\Omega)^N$  there exists a 
unique $\theta \in \CW^1_q(\Omega)$ which satisfies the  
variational equation:
\begin{equation}\label{eq:var1}
(\nabla\theta, \nabla\varphi)_\Omega = (f, \nabla\varphi)_\Omega
\quad\text{for all $\varphi \in\CW^1_{q'}(\Omega)$,}
\end{equation} 
and the estimate: $\|\nabla\theta\|_{L_q(\Omega)} 
\leq C_q\|f\|_{L_q(\Omega)}$ 
for some constant $C_q$ independent of $f$, $\theta$ and $\varphi$. 
\end{dfn}
\begin{remark}\label{rem:3} 
\thetag1~$W^1_q(\Omega) + \CW^1_q(\Omega) = \{p = p_1 + p_2 \mid 
p_1 \in W^1_q(\Omega), \enskip p_2 \in \CW^1_q(\Omega)\}$  
is the space for pressures. 
\vskip0.5pc\noindent
\thetag2~When $\Omega$ is a bounded domain, a half-space, a perturbed
half-space, or a layer domain, 
 weak Dirichlet-Neumann problem is uniquely solvable 
 with  $\CW^1_q(\Omega) 
= \hat W^1_{q, 0}(\Omega)$, while  
 when $\Omega$ is an exterior domain, it is uniquely solvable 
with $\CW^1_q(\Omega)$ being the closure
of $W^1_{q, 0}(\Omega)$ by semi-norm $\|\nabla\cdot\|_{L_q(\Omega)}$.
More examples of domains where the unique solvability of 
weak Dirichlet-Neumann problem holds were given in \cite{S0, S1}.
\end{remark}

To state the compatibility condition for initial data $\bv_0$, we introduce
the solenoidal space $J_q(\Omega)$ defined by 
$J_q(\Omega) = \{\bff \in L_q(\Omega)^N \mid (\bff, \nabla\varphi)_\Omega
= 0 \enskip\text{for any $\varphi \in \CW^1_{q'}(\Omega)$}\}$. 
Since $C^\infty_0(\Omega) \subset \CW^1_{q'}(\Omega)$, 
we see that $\dv \bff = 0$ in $\Omega$ provided that 
$\bff \in J_q(\Omega)$.  But, the opposite direction does not hold
in general. We define $\CD_{q,p}(\Omega)$ by $\CD_{q,p}(\Omega)
=(J_q(\Omega), \CD_q(\Omega))_{1-1/p. p}$ with 
\begin{gather}
\CD_q(\Omega) = \{\bff \in W^2_q(\Omega)^N \mid  
\text{$\bff$ satisfies the compatibility condition:} \nonumber\\
(\bD(\bff)\bn - <\bD(\bff)\bn, \bn>\bn)|_{S} = 0,
\quad \bff|_{\Gamma} = 0\}.\label{comp:2}
\end{gather}
From Steiger \cite{Steiger}, we know that 
$$\CD_{q,p}(\Omega) = 
\begin{cases}
\{\bff \in B^{2(1-1/p)}_{q,p}(\Omega) \cap J_q(\Omega) \mid 
\text{$\bff$ satisfies \eqref{comp:2}}\} 
&\quad\text{when $2(1-1/p) > 1+1/q$}, \\
\{ \bff \in J_q(\Omega) \cap B^{2(1-1/p)}_{q,p}(\Omega) 
\mid \bff|_\Gamma = 0\} &\quad\text{when 
$1/q < 2(1-1/p) < 1+1/q$}, \\
B^{2(1-1/p)}_{q,p}(\Omega) \cap J_q(\Omega)
&\quad\text{when 
$2(1-1/p) < 1/q$}.
\end{cases}
$$
The following theorem is concerned with local in time unique 
existence theorem for \eqref{p1}.
\begin{thm}\label{main:loc} Let $2 < p < \infty$,  
$N < q < \infty$ and $R > 0$. Assume that 
$\Omega$ is a uniform $W^{2-1/q}_q$ domain and that 
 weak Dirichlet-Neumann problem is uniquely solvable 
for $\CW^1_q(\Omega)$ and $\CW^1_{q'}(\Omega)$
$(q' = q/(q-1))$. Then, there exists a time
$T > 0$ depending on $R$ such that for any initial data
$\bv_0 \in \CD_{q,p}(\Omega)$ 
with $\|\bv_0\|_{B^{2(1-1/p)}_{q,p}(\Omega)} \leq R$
problem \eqref{p1} admits a unique solution 
$\bu \in L_q((0, T), W^2_q(\Omega)) \cap W^1_p((0, T), L_q(\Omega))$
with some pressure term 
$\theta \in L_p((0, T), W^1_q(\Omega) + \CW^1_q(\Omega))$ 
possessing the estimate:
$$\|\bu\|_{L_p((0, T), W^2_q(\Omega))} 
+ \|\pd_t\bu\|_{L_p((0, T), L_q(\Omega))} \leq M_0R$$
with some positive constant $M_0$ independent of $R$ and $T$. 
\end{thm}
\begin{remark}\label{rem:5}
\thetag1~ Employing the similar argumentation to Str\"omer
\cite{St}, we can prove that there exists a positive number $\sigma > 0$
such that the map: $x = \bX_\bu(\xi, t)$ is diffeomorphism
from $\Omega$ onto $\Omega_t$, $S$ onto $S_t$ and $\Gamma$ onto
$\Gamma$ for any $t \in (0, T)$ provided that 
\begin{equation}\label{small:1}
\int^T_0\|\nabla \bu(\cdot, t)\|_{L_\infty(\Omega)}\,dt \leq 
\sigma,
\end{equation}
 so that from Theorem \ref{main:loc} 
$\bv(x, t) = \bu(\bX_\bu^{-1}(x, t), t)$ 
solves the original free boundary problem \eqref{1.1}-\eqref{1.2}
for small $T>0$ with some pressure term $\pi$, 
where $\bX^{-1}_\bu(x, t)$
denotes the inverse map of the correspondence: $x = \bX_\bu(\xi, t)$.
\vskip0.5pc\noindent
\thetag2~ It is easy to extend Theorem
\ref{main:loc} to  the equation:
\begin{equation}\label{right:1}
\pd_t\bv + (\bv\cdot\nabla)\bv - \DV\bT(\bv, \theta) = \bff, \quad
\dv \bv = 0
\end{equation}
instead of \eqref{1.1} under similar assumption on $\bff$
to Solonnikov \cite{Sol1} and Shibata and Shimizu
\cite{SS1}. But,  we only consider
the case $\bff = 0$ in this paper for simplicity. 
\end{remark}

Our  global in time unique existence theorem is obtained under the assumption 
that $\Omega$ is a bounded domain and the key issue is the orthogonality of 
the rigid motion.  We introduce the 
rigid space $\CR_d$ defined by 
\begin{equation}\label{rigid:1}
\CR_d = \{Ax+b \mid A: N\times N \enskip \text{anti-symmetric matrix},
 b \in \BR^N\}.
\end{equation}
We know that $\bu$ satisfies $\bD(\bu) =0$ if and only if $\bu \in \CR_d$
(cf. \cite{DL}).  Let $\{\bp_\ell\}_{\ell=1}$ 
be the orthogonal bases of $\CR_d$, that is
$\bp_\ell \in \CR_d$ $(\ell=1, \ldots, M)$ and 
\begin{equation}\label{rigid:2} (\bp_\ell, \bp_m)_\Omega = \delta_{\ell m}
\quad(\ell, m=1, \ldots, M),
\end{equation}
where $\delta_{\ell m}$ is the Kronecker delta symbol such that 
$\delta_{\ell\ell} = 1$ and $\delta_{\ell m} = 0$ with $\ell\not=m$, $M$ 
the dimension of $\CR_d$ and $(\cdot, \cdot)_\Omega$ the $L_2$ inner-product 
on $\Omega$. The following theorm is 
our global in time unique existence result.
\begin{thm}\label{main:global} Let $2 < p < \infty$ and $N < q < \infty$.
Assume that $\Omega$ is a bounded domain, 
$S$ and $\Gamma$ are  $W^{2-1/q}_q$ compact hypersurfaces
and that $S \not=\emptyset$.  Then, there exist  numbers  
$\epsilon > 0$ and $\gamma > 0$ 
such that for any initial data $\bv_0 \in \CD_{q,p}
(\Omega)$ with $\|\bv_0\|_{B^{2(1-1/p)}_{q,p}(\Omega)} \leq \epsilon$
that satisfies, in
addition, the orthogonality condition:
\begin{equation}\label{orth:1}
(\bv_0, \bp_\ell)_\Omega= 0 \quad(\ell=1, \ldots, M)
\quad\text{when $\Gamma = \emptyset$,}
\end{equation}
 problem \eqref{p1} with $T = \infty$
 admits a unique
solution $\bu \in L_p((0, \infty), W^2_q(\Omega)) 
\cap W^1_p((0, \infty), L_q(\Omega))$ 
possessing the estimate:
$$\|e^{\gamma t}\bu\|_{L_p((0, \infty), W^2_q(\Omega))}
+ \|e^{\gamma t}\pd_t\bu\|_{L_p((0, \infty), L_q(\Omega))}
\leq C\epsilon. $$
for some positive constant $C$ independent of $\epsilon$. 
\end{thm}
\section{A proof of a local in time unique existence theorem}

In this section, we prove Theorem \ref{main:loc}. For this purpose, first we 
state our maximal $L_p$-$L_q$ regularity theorem obtained by Shibata 
\cite{S1} for the linearized system \eqref{lp}. 
To state our maximal regularity result for \eqref{lp}, we introduce 
some symbols.  For any Banach space $X$ with norm $\|\cdot\|_X$,  
integer $m \geq 0$ and $\gamma_0 > 0$, we set
\begin{align*}
W^m_{p, \gamma_0}(I, X) & = \{f : I \to X \mid 
e^{-\gamma_0 t}f(t) \in W^m_p(\BR_+, X)\}\quad(I = \BR_+, \BR),\\
W^m_{p, 0, \gamma_0}(\BR, X) & = \{f : \BR \to X \mid 
e^{-\gamma_0 t}f(t) \in W^m_p(\BR, X), \enskip 
f(t) = 0\enskip\text{for $t < 0$}\},
\end{align*}
where $\BR_+ = (0, \infty)$. We set $W^0_{p, \gamma_0} = L_{p, \gamma_0}$
and $W^0_{p, 0, \gamma_0} = L_{p,0,\gamma_0}$. Let $\CL$ and 
$\CL^{-1}$ be the Laplace transform and its inverse transform defined by 
$$\CL[f](\lambda) = \int^\infty_{-\infty}e^{-\lambda t}f(t)\,dt,\quad
\CL^{-1}[g](t) = \frac{1}{2\pi}\int^\infty_{-\infty}
e^{\lambda t}g(\gamma + i\tau)\,d\tau$$
with $\lambda = \gamma + i\tau \in \BC$. For any real number $s \geq 0$, 
let $H^s_{p, \gamma_0}(\BR, X)$ be the Bessel potential space of order
$s$ defined by 
$$H^s_{p,\gamma_0}(\BR, X) = \{f \in L_{p, \gamma_0}(\BR, X) \mid
e^{-\gamma t}\Lambda^s_\gamma f \in L_p(\BR, X)
\enskip\text{for any $\gamma \geq \gamma_0$}\}
$$
with $[\Lambda^s_\gamma f](t) 
= \CL^{-1}[\lambda^s\CL[f](\lambda)](t)$. We set 
$H^s_{p,0,\gamma_0}(\BR, X)= \{f \in H^s_{p, \gamma_0}(\BR, X) \mid
f(t) = 0 \enskip\text{for $t < 0$}\}$. 
By using the $\CR$ bounded solution opeator 
$\CR(\lambda)$ introduced in Sect.1,  Shibata \cite{S1}
proved the following maximal $L_p$-$L_q$ result for problem 
\eqref{lp}.
\begin{thm}\label{thm:linear} Let $1 < p, q < \infty$, $N < r < \infty$
and $\max(q, q') \leq r$ $(q' = q/(q-1))$.  Assume that $\Omega$ is a 
uniform $W^{2-1/r}_r$ domain and that the weak Dirichlet-Neumann 
problem is uniquely solvable for $\CW^1_q(\Omega)$ and $\CW^1_{q'}(\Omega)$
$(q' = q/(q-1))$. 
Then, there exists a positive number $\gamma_0$ such that
for any initial data $\bu_0 \in \CD_{q,p}(\Omega)$
and right members $\bff$, $g$, $\bg$ and $\bh$ with 
\begin{alignat*}2
\bff &\in L_{p, 0, \gamma_0}(\BR, L_q(\Omega)^N), &\enskip 
g &\in L_{p, 0, \gamma_0}(\BR, W^1_q(\Omega)) \cap 
H^{1/2}_{p, 0, \gamma_0}(\BR, L_q(\Omega)), \\
\bg &\in W^1_{p, 0, \gamma_0}(\BR, L_q(\Omega)^N), &\enskip 
\bh &\in L_{p, 0, \gamma_0}(\BR, W^1_q(\Omega)^N) \cap 
H^{1/2}_{p, 0, \gamma_0}(\BR, L_q(\Omega)^N),
\end{alignat*}
problem \eqref{lp}  admits  a unique solution 
$\bu \in L_{p, \gamma_0}(\BR_+, W^2_q(\Omega)) \cap 
W^1_{p, \gamma_0}(\BR_+, L_q(\Omega)^N)$ with some pressure
term $\theta \in L_{p, \gamma_0}(\BR_+, W^1_q(\Omega) + \CW^1_q(\Omega))$ 
possessing the estimate:
\begin{align*}
&\|e^{-\gamma t}\pd_t\bu\|_{L_p(\BR_+, L_q(\Omega))}
+ \|e^{-\gamma t}\bu\|_{L_p(\BR_+, W^2_q(\Omega))}
\leq C\{\|\bu_0\|_{B^{2(1-1/p)}_{q,p}(\Omega))} \\
&\quad + 
\|e^{-\gamma t}(\bff, \pd_t\bg)\|_{L_p(\BR, L_q(\Omega))}
+ \|e^{-\gamma t}(g, \bh)\|_{L_p(\BR, W^1_q(\Omega))}
+ \|e^{-\gamma t}\Lambda^{1/2}_\gamma(g, \bh) 
\|_{L_p(\BR, L_q(\Omega))}\}
\end{align*}
for any $\gamma \geq \gamma_0$ with some constant $C$ independent of 
$\gamma \geq \gamma_0$. 
\end{thm}
To prove Theorem \ref{main:loc}, we use the maximal $L_p$-$L_q$ regularity
theorem for problem \eqref{lp} in a finite time interval, which is 
derived from Theorem \ref{thm:linear}.  But,
 we have to replace the nonlocal operator $\Lambda^{1/2}_\gamma$ 
with value in $L_q(\Omega)$ 
by the local operator $\pd_t$ with value in $W^{-1}_q(\Omega)$. For this 
purpose, first of all, we introduce the extension map
$\iota : L_{1, {\rm loc}}(\Omega) \to L_{1, {\rm loc}}(\BR^N)$  
having the following properties:
\begin{itemize}
\item[\thetag{e-1}]~ For any $1 < q < \infty$ and $f \in W^1_q(\Omega)$, 
$\iota f \in W^1_q(\BR^N)$, $\iota f = f$ in $\Omega$ 
and $\|\iota f\|_{W^i_q(\BR^N)} \leq C_q\|f\|_{W^i_q(\Omega)}$ 
for $i=0, 1$ with some constant $C_q$ depending on $q$, $r$
 and $\Omega$. 
\item[\thetag{e-2}]~ For any $1 < q < \infty$ and $f \in W^1_q(\Omega)$,
$\|(1-\Delta)^{-1/2}\iota(\nabla f)\|_{L_q(\BR^N)} \leq C_q\|f\|_{L_q(\Omega)}$
with some constant $C_q$ depending on $q$, $r$ and $\Omega$.
\end{itemize}
Here, $(1-\Delta)^{-1/2}$ is the operator defined by 
$(1-\Delta)^{-1/2}f = \CF^{-1}[(1+|\xi|^2)^{1/4}\CF[f]]$ with the help of 
Fourier transform
$\CF$ and Fourier inverse transform $\CF^{-1}$ which are defined by 
$$\CF[f](\xi) = \int_{\BR^N} e^{-ix\cdot\xi}f(x)\,dx,
\quad \CF^{-1}[g(\xi)] = \frac{1}{(2\pi)^N}\int_{\BR^N} 
e^{ix\cdot\xi}g(\xi)\,d\xi.$$
In the following, such extension map $\iota$ is fixed. We define 
$W^{-1}_q(\Omega)$ by 
$$W^{-1}_q(\Omega) = \{f \in L_{1, {\rm loc}}(\Omega) \mid 
(1- \Delta)^{-1/2}\iota f \in L_q(\Omega)\}.$$
As is proved in the appendix below, we have
\begin{gather}
W^1_{p, 0, \gamma_0}(\BR, W^{-1}_q(\Omega)) 
\cap L_{p, 0, \gamma_0}(\BR, W^1_q(\Omega)) 
\subset H^{1/2}_{p, 0, \gamma_0}(\BR, L_q(\Omega)),
\label{emb:1} \\
\|e^{-\gamma t}\Lambda^{1/2}_\gamma f\|_{L_p(\BR, L_q(\Omega))}
\leq C\{\|e^{-\gamma t}\pd_t[(1-\Delta)^{-1/2}(\iota f)]
\|_{L_p(\BR, L_q(\BR^N))} + \|e^{-\gamma t}f
\|_{L_p(\BR, W^1_q(\Omega))}\} \label{emb:2}
\end{gather}
for any $\gamma\geq \gamma_0$. Combining Theorem \ref{thm:linear}
with \eqref{emb:1}, we have the following theorem.
\begin{thm}\label{thm:linear*} 
Let $1 < p, q < \infty$, $N < r < \infty$
and $\max(q, q') \leq r$ $(q' = q/(q-1))$. 
Let $T$ be any positice number.  Assume that $\Omega$ is a 
uniform $W^{2-1/r}_r$ domain and that  weak Dirichlet-Neumann 
problem is uniquely solvable for $\CW^1_q(\Omega)$ and $\CW^1_{q'}(\Omega)$
$(q' = q/(q-1))$. 
Then, there exists a positive number $\gamma_0$ such that
for any initial data $\bu_0 \in \CD_{q,p}(\Omega)$ and
any right members $\bff$, $g$, $\bg$ and $\bh$ with 
\begin{alignat*}2
\bff &\in L_p((0, T), L_q(\Omega)^N)), &\enskip
g &\in L_p((0, T), W^1_q(\Omega)) \cap W^1_p((0, T), W^{-1}_q(\Omega)), \\
\bg &\in W^1_p((0, T), L_q(\Omega)^N), &\enskip
\bh &\in L_p((0, T), W^1_q(\Omega)^N) \cap W^1_p((0, T), W^{-1}_q(\Omega)^N),
\end{alignat*}
satisfying the conditions: $g|_{t=0} =0$, $\bg|_{t=0} = 0$  
and $\bh|_{t=0} = 0$, problem \eqref{lp} admits a 
unique solution $\bu \in L_p((0, T), W^2_q(\Omega)^N) 
\cap W^1_p((0, T), L_p(\Omega)^N)$ with pressure term 
$\theta \in L_p((0, T), W^1_q(\Omega) + \CW^1_q(\Omega))$ 
possessing the estimate:
\begin{multline*}
\|\bu\|_{L_p((0, t), W^2_q(\Omega))}
+ \|\pd_t\bu\|_{L_p((0, t), L_q(\Omega))}
\leq Ce^{\gamma t}\{
 \|\bu_0\|_{B^{2(1-1/p)}_{q,p}(\Omega)} \\
+ \|(\bff, \pd_t\bg)\|_{L_p((0, t), L_q(\Omega))}
+ \|(g, \bh)\|_{L_p((0, t), W^1_q(\Omega))}
+ \|\pd_t[(1-\Delta)^{-1/2}(\iota g, \iota \bh)]\|_{L_p((0, t),L_q(\BR^N))}\}
\end{multline*}
for any $t \in (0, T]$ and $\gamma \geq \gamma_0$ with some constant 
$C$ independent of $\gamma \geq \gamma_0$ and $t \in (0, T]$.
\end{thm}
\begin{proof}
Let $t$ be any number with $0 < t \leq T$. Given $f(\cdot, s)$ defined for 
$s \geq 0$,
$f_0(\cdot, s)$ denotes the zero extension of $f$ to $s < 0$, that is 
$f_0(\cdot, s) = f(\cdot, s)$ for $s \geq 0$ and $f_0(\cdot, s) = 0$
for $s < 0$.  Let $E_tf$ be the extension of $f$ defined by 
\begin{equation}\label{ext:1}
E_tf = \begin{cases} f_0(\cdot, s) &\quad\text{for $s \leq t$}, \\
f_0(\cdot, 2t-s) &\quad\text{for $s \geq t$}.
\end{cases}
\end{equation}
Note that $E_tf$ vanishes for $s \not\in [0, 2t]$.  Moreover, if 
$f|_{s=0}$, then 
\begin{equation}\label{ext:2}
\pd_sE_tf = \begin{cases} 
\pd_sf(\cdot, s) &\quad\text{for $s \leq t$}, \\
-(\pd_sf)(\cdot, 2t-s) &\quad\text{for $s \geq t$},
\\
0 &\quad\text{for $s \not\in [0, 2t]$}.
\end{cases}
\end{equation}
Let $\bu^t = \bv(\cdot, s)$ and $\theta^t = \kappa(\cdot, s)$ 
be solutions to the equations:
\begin{gather}
\pd_s\bv - \DV\bT(\bv, \kappa) = E_t\bff, \quad
\dv \bv = E_tg = \dv( E_t\bg)
 \quad  \text{in $\Omega\times(0, \infty)$}, 
\nonumber\\
\bT(\bv, \kappa)\tilde\bn|_{S}  = E_t\bh|_{S}, 
\quad \bv|_{\Gamma} = 0, 
\quad
\bv|_{t=0}  = \bu_0 \quad\text{in $\Omega$}.
\label{ext:3}
\end{gather}
Since $E_{t_1}f = E_{t_2}f$ for $0 < t_1 , t_2 \leq T$, by the uniqueness 
of solutions yields that $\bu^{t_1}(\cdot, s) = \bu^{t_2}(\cdot, s)$
for $s \in [0, t_1]$ with $0 < t_1 < t_2 \leq T$.  By Theorem \ref{thm:linear}
\begin{multline}
\|e^{-\gamma s}\pd_s\bu^t\|_{L_p(\BR, L_q(\Omega))}
+ \|e^{-\gamma s}\bu^t\|_{L_p(\BR, W^2_q(\Omega))} 
\leq C\{\|\bu_0\|_{B^{2(1-1/p)}_{q,p}(\Omega)}
+ \|e^{-\gamma s}(E_t\bff, \pd_s(E_t\bg))\|_{L_p(\BR, L_q(\Omega))}\\
+ \|e^{-\gamma s}(E_tg, E_t\bh)\|_{L_p(\BR, W^1_q(\Omega))}
+ \|e^{-\gamma s}\Lambda^{1/2}_\gamma(E_t g, E_t\bh)
\|_{L_p(\BR, L_q(\Omega))}\}. \label{ext:4}
\end{multline}
Noting \eqref{ext:2}, we see easily that 
\begin{align}
\|e^{-\gamma s}(E_t\bff, \pd_s(E_t\bg))\|_{L_p(\BR, L_q(\Omega))} 
&\leq 2\|e^{-\gamma s}(\bff, \pd_s\bg)\|_{L_p((0, t), L_q(\Omega))},
\nonumber\\
\|e^{-\gamma s}(E_tg, E_t\bh)\|_{L_p(\BR, W^1_q(\Omega))} 
&\leq 2\|e^{-\gamma s}(g, \bh)\|_{L_p((0, t), W^1_q(\Omega))}.
\label{ext:5}
\end{align}
Moreover, by \eqref{emb:2} and \eqref{ext:2}, we have 
\begin{align}
&\|e^{-\gamma s}\Lambda^{1/2}_\gamma(E_tg, E_t\bh)\|_{L_p(\BR, L_q(\Omega))}
\nonumber\\
&\quad
\leq C\{\|e^{-\gamma s}\pd_s[(1-\Delta)^{-1/2}
(\iota g, \iota \bh)]\|_{L_p((0, t), L_q(\BR^N))}
+ \|e^{-\gamma s}(g, \bh)\|_{L_p((0, t), W^1_q(\Omega))}\}.
\label{ext:6}
\end{align}
Setting $\bu = \bu^T$ and $\theta = \theta^T$, noting that 
$\bu(\cdot, s) = \bu^t(\cdot, s)$ for $0 < s < t$ and 
combining \eqref{ext:4}, \eqref{ext:5} and \eqref{ext:6},
we complete the proof of Theorem \ref{thm:linear*}.
\end{proof}
{\bf A Proof of Theorem \ref{main:loc}}~ 
In the following, we assume that $2 < p < \infty$ and 
$N < q < \infty$, that $\Omega$ is a uniform $W^{2-1/q}_q$ domain
in $\BR^N$ ($N \geq 2$),
and that  weak Dirichlet-Neumann problem is uniquely solvable 
for $\CW^1_q(\Omega)$ and $\CW^1_{q'}(\Omega)$
$(q' = q/(q-1))$.  By Sobolev's imbedding theorem
we have 
\begin{equation}\label{5.1}
W^1_q(\Omega) \subset L_\infty(\Omega), \quad
\|\prod_{j=1}^m f_j\|_{W^1_q(\Omega)} \leq 
C\prod_{j=1}^m\|f_j\|_{W^1_q(\Omega)}.
\end{equation}
Let $T$ and $L$ be any positive numbers and we define a space $\CI_{L, T}$
by
\begin{equation}\label{5.9}\begin{split}
\CI_{L,T} =\{&\bv \in L_p((0, T), W^2_q(\Omega)) \cap W^1_p((0, T), 
L_q(\Omega)) \mid  \bv|_{t=0} = \bv_0 \quad\text{in $\Omega$}, \quad
\BI_\bv(0,T) \leq L\},
\end{split}\end{equation}
where we have set
$\BI_\bv(0,T) = \|\bv\|_{L_p((0, T), W^2_q(\Omega))}
+ \|\pd_t\bv\|_{L_p((0, T), L_q(\Omega))}$. 
Given $\bw \in \CI_{L, T}$, let $\bv$ and $\omega$ be solutions to
problem:
\begin{equation}\label{5.8}\begin{split}
\pd_t\bv - \DV\bT(\bv, \omega) = \bF(\bw), \enskip 
\dv\bv = G(\bw) = \dv \bG(\bw)
\quad&\text{in $\Omega\times
(0, T)$}, \\
\bT(\bv, \omega)\tilde\bn|_{S} = \bH(\bw)\tilde\bn|_{S},
\quad\bv|_{\Gamma} = 0, \quad
\bv|_{t=0} = \bv_0\quad&\text{in $\Omega$}.
\end{split}
\end{equation}
First, we estimate the right-hand sides of \eqref{5.8}. 
By \eqref{5.1}  and H\"older's inequality we have
\begin{equation}\label{eq:2.1}\begin{split}
\sup_{t \in (0, T)}\|\int^t_0\nabla\bw(\cdot, s)\,ds
\|_{L_\infty(\Omega)}  \leq M_1T^{1/{p'}}L, \quad
\sup_{t \in (0, T)}\|\int^t_0\nabla\bw(\cdot, s)\,ds
\|_{W^1_q(\Omega)}  \leq CT^{1/{p'}}L.
\end{split}\end{equation}
with $p' = p/(p-1)$. 
Here and in the following, $C$ denotes a generic constant independent of 
$T$ and $R$ and we use the letters $M_i$ 
to denote some special constants independent of $T$ and $L$. 
The value of $C$ may change from line to line. 
To treat nonlinear functions with respect to $\int^t_0 \nabla\bw(\cdot, s)\,ds$, we choose 
$T$ so small that $M_1T^{1/{p'}}L \leq 1$ in 
\eqref{eq:2.1}, so that
\begin{equation}\label{eq:2.2}
\sup_{t \in (0, T)} \|\int^t_0\nabla\bw(\cdot, s)\,ds
\|_{L_\infty(\Omega)} \leq 1. 
\end{equation}
By \eqref{eq:2.1}, \eqref{eq:2.2}, \eqref{5.1} and \eqref{null:1}, we have
\begin{equation}\label{eq:2.3}
\sup_{t \in (0, T)}\|\bV_i(\int^t_0\nabla\bw(\cdot, t)\,ds)\|_{W^1_q(\Omega)}
\leq CT^{1/{p'}}L, \quad 
\sup_{t \in (0, T)}\|\nabla\bW(\int^t_0\nabla\bw(\cdot, t)\,ds)\|_{L_q(\Omega)}
\leq CT^{1/{p'}}L
\end{equation}
where $i=1,2,4,5$ and $6$, and $\bW = \bW(\bK)$ is any matrix of polynomials
with respect to $\bK$. By \eqref{nonlinear:1},  
\eqref{ext-normal}, \eqref{5.1}, \eqref{eq:2.1}, \eqref{eq:2.2} and 
\eqref{eq:2.3}, we have
\begin{equation}\label{eq:2.6*}
\|\bF(\bw)\|_{L_p((0, T), L_q(\Omega))}
\leq CL^2T^{1/{p'}}, \quad 
\|(G(\bw), \bH(\bw)\tilde\bn)\|_{L_p((0, T), W^1_q(\Omega))}
\leq CL^2T^{1/{p'}}. 
\end{equation}

To obtain 
\begin{equation}\label{eq:2.4}
\sup_{t \in (0, T)}\|\bw(\cdot, t)\|_{B^{2(1-1/p)}_{q,p}(\Omega)}
\leq C(\BI_\bw(0,T) + e^{\gamma T}
\|\bv_0\|_{B^{2(1-1/p)}_{q,p}(\Omega)}),
\end{equation}
we use 
 the embedding relation:
\begin{equation}\label{SS:2.18}
L_p((0, \infty), X_1) \cap W^1_p((0, \infty), X_0) 
\subset BUC(J, [X_0, X_1]_{1-1/p,p})
\end{equation}
for any two Banach spaces $X_0$ and $X_1$ such that $X_1$
is dense in $X_0$ and $1 < p < \infty$
 (cf. \cite{Amann2}). In fact, let $E_t$ be the extension operator
defined in the proof of Theorem \ref{thm:linear*} and 
let $\bZ$ and $\Pi$ be solutions to problem:
\begin{gather}
\pd_t\bZ - \DV\bT(\bZ, \Pi) = 0, \quad\dv\bZ = 0 \quad
\text{in $\Omega\times(0, \infty)$}, \nonumber\\
\bT(\bZ, \Pi)\tilde\bn|_{S} = 0, \quad\bZ|_{\Gamma} = 0, \quad
\bZ|_{t=0} = \bv_0 \quad\text{in $\Omega$}.\label{5.3}
\end{gather}
By Theorem \ref{thm:linear} \thetag1, we know the unique existence of
$(\bZ, \Pi)$ possessing the estimate:
\begin{equation}\label{5.6}
\|e^{-\gamma t}\pd_t\bZ\|_{L_p(\BR_+,L_q(\Omega))} + \|e^{-\gamma t}
\bZ\|_{L_p(\BR_+, W^2_q(\Omega))}
\leq C\|\bv_0\|_{B^{2(1-1/p)}_{q,p}(\Omega)} \quad(\gamma \geq \gamma_0)
\end{equation}
for some constants $\gamma_0$ and $C$, where $C$ is independent of $\gamma
\geq \gamma_0$. We choose $\gamma$ so large and fix it in the following. 
Set $\bz = \bw - \bZ$. Since $\bz|_{t=0} = 0$, 
by \eqref{ext:1} and \eqref{ext:2} we have 
$$\BI_{E_T\bz}(0,\infty) \leq C\BI_\bz(0,T)
\leq C(\BI_\bw(0,T) + e^{\gamma T}\BI_{e^{-\gamma t}\bZ}(0,\infty)).
$$ 
Thus, noting that $\bw = \bZ + E_T\bz$ for 
$t \in (0, T)$ and using \eqref{SS:2.18}, we have
\begin{align*}
&\sup_{t\in(0, T)}\|\bw(\cdot, t)\|_{B^{2(1-1/p)}_{q,p}(\Omega)}
\leq \sup_{t \in (0, \infty)}\| E_T\bz(\cdot, t)\|_{B^{2(1-1/p)}_{q,p}(\Omega)}
+ e^{\gamma T}\sup_{t \in (0, \infty)}\|e^{-\gamma t}\bZ(\cdot, t)
\|_{B^{2(1-1/p)}_{q,p}(\Omega)}\\
&\quad
\leq C(\BI_{E_T\bz}(0,\infty) 
+ e^{\gamma T}\BI_{e^{-\gamma t}\BZ}(0,\infty))
\leq C(\BI_{\bw}(0,T) + e^{\gamma T}\BI_{e^{-\gamma t}\BZ}(0,\infty)),
\end{align*}
which combined with \eqref{5.6} furnishes \eqref{eq:2.4}.  

Since $B^{2(1-1/p)}_{q,p}(\Omega) \subset W^1_q(\Omega)$ as follows 
from the assumption: $2 < p < \infty$, by \eqref{eq:2.4} and \eqref{eq:2.2} 
we have 
\begin{equation}\label{eq:2.5}\begin{split}
\sup_{t \in (0, T)}\|\bw(\cdot, t)\|_{W^1_q(\Omega)} &
\leq C(L+e^{\gamma T}\|\bv_0\|_{B^{2(1-1/p)}_{q,p}(\Omega)}), \\
\sup_{t \in (0, T)}\|\pd_t\bW(\int^t_0\nabla\bw(\cdot, s)\,ds)\|_{L_q(\Omega)}
&\leq C(L+e^{\gamma T}\|\bv_0\|_{B^{2(1-1/p)}_{q,p}(\Omega)}).
\end{split}\end{equation}
Writing $\pd_t\bG(\bw) = \{\pd_t\bV_5(\int^t_0\nabla \bw)\,ds)\}\bw
+ \bV_5(\int^t_0\nabla \bw\,ds)\pd_t\bw$
and using \eqref{null:1}, \eqref{5.1}, \eqref{eq:2.4}, \eqref{eq:2.5} 
and \eqref{eq:2.3}, we have
\begin{equation}\label{eq:2.7}
\|\pd_t\bG(\bw)\|_{L_p((0, T), L_q(\Omega))}
\leq C\{L^2T^{1/{p'}} + 
(L + e^{\gamma T}\|\bv_0\|_{B^{2(1-1/p)}_{q,p}(\Omega)})^2T^{1/p})\}.
\end{equation}


To continue our estimate, we prepare the following lemma.
\begin{lem}\label{lem:est1}
Let $1 < p < \infty$, $N < q, r < \infty$ and let $\Omega$ be
a uniform $W^{2-1/r}_r$ domain. Let $\iota$ be the 
extension map satisfying the properties \thetag{e-1}
and \thetag{e-2}. Then, 
\begin{align*}
&\|\pd_t[(1-\Delta)^{-1/2}\iota((\nabla f)g)]\|_{L_p((0, T), L_q(\BR^N))}
\\
&\quad
\leq C\Bigl\{\Bigl(\int^T_0(\|\pd_tf(\cdot, t)\|_{L_q(\Omega)}
\|g(\cdot, t)\|_{W^1_q(\Omega)})^p\,dt\Bigr)^{1/p}
+ \Bigl(\int^T_0(\|\nabla f(\cdot, t)\|_{L_q(\Omega)}
\|\pd_tg(\cdot, t)\|_{L_q(\Omega)})^p\,dt\Bigr)^{1/p}
\Bigr\}.
\end{align*}
\end{lem}
\begin{proof}
To prove the lemma, we use 
an inequality:
\begin{equation}\label{eq:2.6}
\|(1-\Delta)^{-1/2}\iota(fg)\|_{L_q(\BR^N)}
\leq C\|f\|_{L_q(\Omega)}\|g\|_{L_q(\Omega)}
\end{equation}
provided that $N < q < \infty$, which follows from 
the following observation:   For any $\varphi
\in C^\infty_0(\BR^N)$ by H\"older's inequality and \thetag{e-1}
we have 
\begin{align*}
|((1-\Delta)^{-1/2}\iota(fg), \varphi)_{\BR^N}|
=|(\iota(fg), (1-\Delta)^{-1/2}\varphi)_{\BR^N}|
\leq C\|f\|_{L_q(\Omega)}\|g\|_{L_q(\Omega)}
\|(1-\Delta)^{-1/2}\varphi\|_{L_s(\BR^N)},
\end{align*}
where $s$ is an index such that $2/q+1/s = 1$.  Since 
$N(1/{q'} - 1/s) = N/q < 1$, by Sobolev's imbedding theorem 
we have $\|(1-\Delta)^{-1/2}\varphi\|_{L_s(\BR^N)}
\leq C\|\varphi\|_{L_{q'}(\BR^N)}$, which furnishes \eqref{eq:2.6}. 
Since 
\begin{align*}
\pd_t[(1-\Delta)^{-1/2}\iota((\nabla f)g)]
= (1-\Delta)^{-1/2}\iota[\nabla\{(\pd_tf)g\}] - (1-\Delta)^{-1/2}
\iota[(\pd_tf)(\nabla g)]
+ (1-\Delta)^{-1/2}\iota[(\nabla f)\pd_tg], 
\end{align*}
by \eqref{eq:2.6}, \eqref{5.1} and \thetag{e-2} we have
Lemma \ref{lem:est1}.
\end{proof}
Applying Lemma \ref{lem:est1} to $G(\bw)$ and 
$\bH(\bw)\tilde\bn$ with 
$f = \bw$, $g = \bV_4(\int^t_0\nabla(\bw)\,ds)$ and 
$f = \bw$, $g = \bV_6(\int^t_0\nabla(\bw)\,ds)\tilde\bn$,
respectively, and using \eqref{eq:2.3}, 
\eqref{eq:2.4} and \eqref{eq:2.5},
we have
\begin{equation}\label{eq:2.8}\begin{split}
\|\pd_t[(1-\Delta)^{-1/2}&(\iota G(\bw), 
\iota(\bH(\bw)\tilde\bn))]\|_{L_p((0, T), 
L_q(\BR^N))} \\
&\leq C\{L^2T^{1/{p'}} + (L + 
e^{\gamma T}\|\bv_0\|_{B^{2(1-1/p)}_{q,p}(\Omega)})^2T^{1/p}\}.
\end{split}\end{equation}
Thus, applying Theorem \ref{thm:linear*}
to problem \eqref{5.8} and 
using \eqref{eq:2.6*}, \eqref{eq:2.7} and \eqref{eq:2.8}, we have 
\begin{equation}\label{eq:2.9}
\BI_\bv(0,T) \leq M_2\|\bv_0\|_{B^{2(1-1/p)}_{q,p}(\Omega)} + 
M_3 (L^2T^{1/{p'}} 
+ (L + e^{\gamma T}\|\bv_0\|_{B^{2(1-1/p)}_{q,p}(\Omega)})^2T^{1/p}).
\end{equation} 
Let $R$ be a number such that $\|\bv_0\|_{B^{2(1-1/p)}_{q,p}(\Omega)} \leq R$
and set $L = (M_2+1)R$.  Choosing $T > 0$ so small 
that 
$$M_3 (L^2T^{1/{p'}} 
+ (L + e^{\gamma T}R)^2T^{1/p})
\leq 1,
$$
by \eqref{eq:2.9} we have $\BI_\bv(0,T) \leq L$,
so that $\bv \in \CI_{L, T}$. 
If we define a map $\Phi$ by $\Phi(\bw) =
\bv$, then $\Phi$ is a map from 
$\CI_{L,T}$ into itself. 

Next, we show the contractility of the map $\Phi$ on
$\CI_{L,T}$. Let $\bw_i \in \CI_{L,T}$ and set 
$\bv_i = \Phi(\bw_i)$ ($i = 1, 2$).  Setting $\bv = \bv_1 - \bv_2$,
we have
\begin{gather}
\pd_t\bv - \DV\bT(\bv, \omega) = \bff(\bw_1, \bw_2),
\enskip \dv\bv = g(\bw_1, \bw_2) = \dv \bg(\bw_1, \bw_2)
\quad\text{in $\Omega\times(0, T)$}, \nonumber\\
\bT(\bv, \omega)\tilde\bn|_S = \bh(\bw_1, \bw_2)|_S, 
\enskip \bv|_\Gamma = 0, \enskip 
\bv|_{t=0} \enskip\text{in $\Omega$}
\label{eq:2.10}
\end{gather}
with some pressure term $\omega$, where we have set 
\begin{alignat*}2
\bff(\bw_1, \bw_2) & = 
\bF(\bw_1) - \bF(\bw_2), &\quad g(\bw_1, \bw_2) &= 
G(\bw_1) - G(\bw_2), \\
\bg(\bw_1, \bw_2) & = \bG(\bw_1)- \bG(\bw_2),
&\quad 
\bh(\bw_1, \bw_2) &= (\bH(\bw_1) - \bH(\bw_2))\tilde\bn
\end{alignat*}
By Theorem \ref{thm:linear*}, we have 
\begin{equation}\label{eq:2.11}
\BI_{\bv_1- \bv_2}(0,T) \leq M_4\,\BJ(\bw_1, \bw_2)(T)
\end{equation}
for some constant $M_4$ independent of $T$ and $R$ with 
\begin{align*}
\BJ(\bw_1, \bw_2) &= 
\|(\bff(\bw_1, \bw_2), \pd_t\bg(\bw_1, \bw_2))\|_{L_p(0, T), L_q(\Omega))}
+ \|(g(\bw_1, \bw_2), \bh(\bw_1, \bw_2))\|_{L_p((0, T), W^1_q(\Omega))}\\
&\quad + \|\pd_t[(1-\Delta)^{-1/2}(\iota g(\bw_1, \bw_2), \iota
\bh(\bw_1, \bw_2))]\|_{L_p((0, T), L_q(\BR^N))}.
\end{align*}
We estimate each terms in the right-hand side of \eqref{eq:2.10}.
Recalling that  $\|\bv_0\|_{B^{2(1-1/p)}_{q,p}(\Omega)} \leq R$
and $L = (M_2+1)R$, by \eqref{eq:2.1}, \eqref{eq:2.2}, \eqref{eq:2.3}, 
\eqref{eq:2.4} and \eqref{eq:2.5} we have 
\begin{alignat*}2
&\sup_{t\in (0, T)} \|\int^t_0 \nabla\bw_i(\cdot, s)\,ds\|_{L_\infty(\Omega)}
\leq 1, &\quad 
&\sup_{t \in (0, T)}\|\int^t_0\nabla\bw_i(\cdot, s)\,ds\|_{W^1_q(\Omega)}
\leq CRT^{1/{p'}}, \\
&\sup_{t \in (0, T)}\|\bV_j(\int^t_0\nabla\bw_i(\cdot, s)\,ds)
\|_{W^1_q(\Omega)}
\leq CRT^{1/{p'}}, &\quad  
&\sup_{t \in (0, T)}\|\nabla\bW(\int^t_0\nabla\bw_i(\cdot, s)\,ds)
\|_{W^1_q(\Omega)}
\leq CRT^{1/{p'}}, \\
&\sup_{t \in (0, T)}\|\bw_i(\cdot, t)\|_{W^1_q(\Omega)}
\leq CRe^{\gamma T}, &\quad 
&\sup_{t \in (0, T)}\|\pd_t\bW(\int^t_0\nabla \bw_i(\cdot, s)\,ds)
\|_{L_q(\Omega)} \leq CRe^{\gamma T},
\end{alignat*}
where $i = 1, 2$ and $j = 1,2,4,5,6$.
Thus, we have  
\begin{align*}
&\sup_{t \in (0, T)} \|\int^t_0\nabla\bw_1(\cdot, s)\,ds
- \int^t_0\nabla\bw_2(\cdot, s)\,ds\|_{W^1_q(\Omega)}
\leq CT^{1/{p'}}\BI_{\bw_1-\bw_2}(0,T), \\
&\sup_{t \in (0, T)} \|\bW(\int^t_0\nabla\bw_1(\cdot, s)\,ds)
- \bW(\int^t_0\nabla\bw_2(\cdot, s)\,ds)\|_{W^1_q(\Omega)}
\leq C(RT^{1/{p'}} + 1)T^{1/{p'}}\BI_{\bw_1-\bw_2}(0,T).
\end{align*}
Since $(\bw_1-\bw_2)|_{t=0} = 0$, employing the similar argumentation
to that in the proof of \eqref{eq:2.4}, we have 
\begin{align*} 
&\sup_{t \in (0, T)} \|\bw_1(\cdot, t) - \bw_2(\cdot, t)\|_{W^1_q(\Omega)}
\leq C\,\BI_{\bw_1-\bw_2}(0,T), \\
&\sup_{t \in (0, T)} \|\pd_t\{\bW(\int^t_0\nabla\bw_1(\cdot, s)\,ds) 
- \bW(\int^t_0\nabla\bw_2(\cdot, s)\,ds)\}\|_{W^1_q(\Omega)}
\leq C(1 + Re^{\gamma T}T^{1/{p'}})\,\BI_{\bw_1-\bw_2}(0,T).
\end{align*}
Using above estimates and Lemma \ref{lem:est1}, we have 
\begin{equation}\label{eq:2.13}
\BJ(\bw_1, \bw_2)(T) \leq M_5\,\CC(R,T)\,\BI_{\bw_1 - \bw_2}(0,T)
\end{equation}
for some constant $M_5$ independent of $R$ and $T$ with 
\begin{align*}
\CC(R, T) &= 
RT^{1/{p'}} + (RT^{1/{p'}})^2
+ (RT^{1/{p'}})^3 + RT^{1/p} \\
& + e^{\gamma T}
\{(RT^{1/p})(RT^{1/{p'}}) + 
RT^{1/p} + RT^{1/{p'}} 
+ (RT^{1/{p'}})^2\} 
 + e^{2\gamma T}(RT^{1/{p'}})(RT^{1/p}).
\end{align*}
Combining \eqref{eq:2.13} with \eqref{eq:2.11} furnishes that 
\begin{equation}\label{eq:2.12}
\BI_{\Phi(\bw_1) - \Phi(\bw_2)}(0,T) 
\leq M_4M_5\,\CC(R,T)\BI_{\bw_1-\bw_2}(0,T).
\end{equation}
Choosing $T$ smaller in such a way that $M_4M_5\CC(R, T) \leq 1/2$, 
we have $\Phi$ is a contraction map on 
$\CI_{L, T}$.  Thus,  the Banach fixed point theorem
tells us that $\Phi$ has a unique fixed point $\bu$ in
$\CI_{L, T}$ satisfying the equations \eqref{p1}.

Finally, we prove the uniqueness.  Given two $\bv_i \in \CI_{L,T}$ $(i=1,2)$
both of which satisfy the equations \eqref{p1} with the same 
initial data $\bv_0 \in B^{2(1-1/p}_{q,p}(\Omega)$, employing the same
argument as in proving \eqref{eq:2.12} and replacing $\bw_i$ by $\bv_i$,
we have $\BI_{\bv_1-\bv_2}(0,T) \leq M_4M_5\,\CC(R, T)\BI_{\bv_1-\bv_2}(0,T)$.
Since  $T$ has been chosen in such a way that 
$M_4M_5\CC(R, T) \leq 1/2$, we have 
$\BI_{\bv_1-\bv_2}(0,T) \leq \frac12\,\BI_{\bv_1-\bv_2}(0,T)$,
which implies that $\bv_1 = \bv_2$. 
This completes the proof of Theorem \ref{main:loc}.
\section{Some decay properties of solutions to problem \eqref{lp}}
In this section, we discuss exponential stability of solutions to problem
\eqref{lp} assuming that $\Omega$ is bounded in addition. 
Let $\CR(\lambda)$ be the $\CR$ bounded solution operator for  
problem \eqref{gs} introduced in Sect. 1.
If we consider the time shifted equation of \eqref{lp}: 
\begin{gather}
\pd_t\bv + \lambda_1\bv - \DV\bT(\bv, \hat \theta) = \bff, \quad
\dv \bv = g = \dv\bg \quad  \text{in $\Omega\times(0, \infty)$}, 
\nonumber\\
\bT(\bv, \hat \theta)\tilde\bn|_{S}  = \bh|_{S}, 
\quad \bv|_{\Gamma} = 0, 
\quad
\bv|_{t=0}  = \bu_0 \quad\text{in $\Omega$},
\label{gs*}
\end{gather} 
a solution $\bv$ is represented by using $\CR(\lambda+\lambda_1)$,
so that  we have 
the following theorem concerning the exponential stability 
of solutions to \eqref{gs*}.
\begin{thm}\label{thm:3.1} Let $1 < p, q < \infty$, $N < r < \infty$
and $\max(q, q') \leq r$ $(q' = q/(q-1))$.  Assume that $\Omega$ is a 
uniform $W^{2-1/r}_r$ domain and that  weak Dirichlet-Neumann 
problem is uniquely solvable for $\CW^1_q(\Omega)$ and $\CW^1_{q'}(\Omega)$
$(q' = q/(q-1))$. 
Then, there exists a $\lambda_1 > \lambda_0$ such that
problem \eqref{gs*} admits a unique solution $(\bv, \hat \theta)$
with 
$$\bv \in L_p(\BR_+, W^2_q(\Omega)) \cap W^1_p(\BR_+, 
L_q(\Omega)), 
\quad \hat\theta \in L_p(\BR_+, W^1_q(\Omega) + \CW^1_q(\Omega))$$
possessing the estimate:
\begin{align*}
&\|e^{\gamma t}\pd_t\bv\|_{L_p(\BR_+, L_q(\Omega))}
+ \|e^{\gamma t}\bv\|_{L_p(\BR_+, W^2_q(\Omega))}\\
&\quad 
 \leq C(\|\bu_0\|_{B^{2(1-1/p)}_{q,p}(\Omega)}+
\|e^{\gamma t}(\bff, \tilde\Lambda^{1/2}_\gamma g, 
\pd_t\bg, \tilde\Lambda^{1/2}_\gamma
\bh)\|_{L_p(\BR, L_q(\Omega))}
+ \|e^{\gamma t}(g, \bh)\|_{L_p(\BR, W^1_q(\Omega))}) 
\end{align*}
for any $\gamma \leq \lambda_0$ with some constants $C$ independent of 
$\gamma \leq \lambda_0$,  provided that $\bu_0 \in \CD_{q,p}(\Omega)$, 
\begin{alignat}3
e^{\gamma t}\bff &\in L_p(\BR, L_q(\Omega)^N), &\enskip 
e^{\gamma t}g &\in L_p(\BR, W^1_q(\Omega)), &\enskip
e^{\gamma t}\tilde\Lambda_\gamma^{1/2}g &
\in L_p(\BR, L_q(\Omega)), \nonumber \\
e^{\gamma t}\pd_t\bg &\in L_p(\BR, L_q(\Omega)^N), &\enskip 
e^{\gamma t}\bh  &\in L_p(\BR, W^1_q(\Omega)^N), &\enskip
e^{\gamma t}\tilde\Lambda_\gamma^{1/2}\bh
&\in L_p(\BR, L_q(\Omega)^N), \label{cond:2}
\end{alignat}
and $\bff$, $g$, $\bg$ and $\bh$ vanish for $t < 0$. Here, we have defined 
$\tilde\Lambda_\gamma^{1/2}f$ by 
\begin{equation}\label{def:2}
\tilde\Lambda_\gamma^{1/2}f = \CL^{-1}_\lambda[(\lambda + \lambda_1)^{1/2}
\CL[f](\lambda)] \quad\text{with $\lambda = -\gamma + i\tau$}.
\end{equation}
\end{thm}

Since the $\CR$ boundedness implies  the 
usual boundedness of operators, 
we also see that for any $\lambda \in \Sigma_{\epsilon, \lambda_0}$, 
$\bff \in L_q(\Omega)^N$, $g \in W^1_q(\Omega)$, $\bg \in L_q(\Omega)^N$ and 
$\bh \in W^1_q(\Omega)^N$, a unique solution  $\bv$ of problem 
\eqref{gs} possesses the generalized resolvent estmate:
\begin{equation}\label{gre}
\|(|\lambda|\bv, |\lambda|^{1/2}\nabla\bv, \nabla^2\bv)\|_{L_q(\Omega)}
\leq C\|(\bff, |\lambda|^{1/2}g, \nabla g, |\lambda|\bg, |\lambda|^{1/2}\bh,
\nabla\bh)\|_{L_q(\Omega)}
\end{equation}
with some constant $C$ depending on $\epsilon$ and $\lambda_0$.
Especially, we see the existence of a continous semigroup
$\{\BT(t)\}_{t\geq 0}$ associated with problem \eqref{lp}, 
which is analytic. 

To prove a global in time unique existence theorem for
\eqref{p1}, we need the exponential stability of solutions to \eqref{lp},
so that from now on, we assume that $\Omega$ is bounded in addition. 
In this case,  weak Dirichlet-Neumann problem is uniquely solvable for any
exponent $q \in (1, \infty)$ with $\CW^1_q(\Omega) = \hat W^1_{q,0}(\Omega)$
and $J_q(\Omega) = \{ \bff \in L_q(\Omega)^N \mid \dv \bff = 0, \enskip  
\bn_{\Gamma}\cdot\bff|_{\Gamma} = 0\}$, where $\bn_{\Gamma}$
is the unit outer normal to $\Gamma$.  
When $\Omega$ is 
bounded, the uniqueness of solutions to problem \eqref{gs} holds 
when $\Gamma\not=\emptyset$ up to $\lambda=0$.  When $\Gamma=\emptyset$,
if we restrict the space of solutions to the quotient space
$W^2_q(\Omega)/\CR_d$, then we also have the uniqueness of solutions to 
\eqref{gs}.  Namely, if $\bu \in W^2_q(\Omega)$ satisfies 
the equations \eqref{gs} with $\bff=0$, $g=0$, $\bg=0$ and $\bh=0$ and 
if $\bu$ satisfies the orthogonal condition: $(\bu, \bp_\ell)_\Omega
= 0$ for $\ell=1, \ldots, M$, then $\bu = 0$ up to $\lambda=0$. 
Moreover, if $\bff \in L_q(\Omega)^N$ and $\bg \in W^1_q(\Omega)^N$
satisfy the condition: $(\bff, \bp_\ell)_\Omega +
<\bh, \bp_\ell>_{S} = 0$, then a solution $\bu$ to problem
\eqref{gs} also satisfies $(\bu, \bp_\ell)_\Omega=0$ whenever $\lambda
\not=0$.  Here, $<f, g>_{S} = \int_\Gamma f(x)g(x)\,d\sigma$,
$d\sigma$ being the surface element of $S$. Using these facts and 
applying a homotopic argument, we see that $\{\BT(t)\}_{t\geq 0}$ 
 is exponentially stable.  Namely, we have the following theorem
which was already proved in Shibata and Shimizu \cite{SS1} in the 
case of $\Gamma=\emptyset$.
\begin{thm}\label{thm:exp-stab} Let $1 <  q < \infty$, $N < r < \infty$
and $\max(q, q') \leq r$ $(q' = q/(q-1))$.  Assume that $\Omega$ is a 
uniform $W^{2-1/r}_r$ domain and that $\Omega$ is bounded in addition. 
Then, there exists a continuous semigroup $\{\BT(t)\}_{t\geq 0}$ 
on $J_q(\Omega)$ associated with problem
\eqref{lp} such that $\bu = \BT(t)\bu_0$ with some pressure 
term $\theta$ solves problem \eqref{lp} 
with $\bff = 0$, $g=0$, $\bg=0$ and $\bh=0$.
Moreover, $\{\BT(t)\}_{t\geq0}$ is analytic and exponentially stable, 
that is 
\begin{equation}\label{exp-stab1}
\|\BT(t)\bu_0\|_{W^\ell_q(\Omega)} \leq C(1+ t^{-\ell/2})e^{-\gamma t}
 \|\bu_0\|_{L_q(\Omega)} \quad\text{for any $t > 0$ and $\ell=0,1,2$}
\end{equation}
with some positive constants $C$ and $\gamma$ 
provided that $\bu_0 \in J_q(\Omega)$ when $\Gamma\not=\emptyset$ and 
$\bu_0 \in J_q(\Omega)$ satisfying the orthogonal condition: 
$(\bu_0, \bp_\ell)_\Omega=0$ for $\ell=1, \ldots, M$ when
$\Gamma=\emptyset$. Here, $W^0_q(\Omega) = L_q(\Omega)$.
\end{thm}
By Theorem \ref{thm:exp-stab}, we have the following 
Corollary which was  proved in Shibata and Shimizu \cite{SS1} in 
the case of $\Gamma=\emptyset$ under the asumption that 
the boundary of $\Omega$ is a $C^{1,1}$ hypersurface.  
\begin{cor}\label{exp-stab2}
Let $1 <  q < \infty$, $N < r < \infty$
and $\max(q, q') \leq r$ $(q' = q/(q-1))$.  Assume that $\Omega$ is a 
uniform $W^{2-1/r}_r$ domain and that $\Omega$ is bounded in addition. 
Then, there exists a  positive constant $\gamma_0$ such that 
problem \eqref{lp} with $\bff=0$, $g=0$, $\bg=0$ and $\bh = 0$ 
admits unique solutions $\bu$ and $\theta$ with 
$$\bu \in L_p(\BR_+, W^2_q(\Omega)^N) \cap W^1_p(\BR_+, L_q(\Omega)^N),
\quad \theta \in L_p(\BR_+, W^1_q(\Omega) + \hat W^1_{q,0}(\Omega))$$
possessing the estimate:
$$
\|e^{\gamma t}\pd_t\bu\|_{L_p(\BR_+,L_q(\Omega))}
+ \|e^{\gamma t}\bu\|_{L_p(\BR_+, W^2_q(\Omega))}
\leq C\|\bu_0\|_{B^{2(1-1/p)}_{q,p}(\Omega)}
\quad\text{for any $\gamma \leq \gamma_0$}  
$$
with some positive constants $C$ independent of $\gamma \leq \gamma_0$ 
provided that $\bu_0 \in \CD_{q,p}(\Omega)$ when $\Gamma\not=\emptyset$ and 
$\bu_0 \in \CD_{q,p}(\Omega)$ and $\bu_0$ satisfies the orthogonal condition: 
$(\bu_0, \bp_\ell)_\Omega=0$ for $\ell=1, \ldots, M$ when
$\Gamma=\emptyset$. 
\end{cor}
Under the preparations mentioned above, we show the following theorem
about the exponential stability of solutions to \eqref{lp}.
\begin{thm}\label{thm:decay} Let $1 < p, q < \infty$, $N < r < \infty$ and 
$\max(q, q') \leq r$ $(q' = q/(q-1))$.  Assume that 
$\Omega$ is a uniform $W^{2-1/r}_r$ domain and that $\Omega$ is bounded
in addition.  Then, there exists a positive constant $\gamma_0$ such that
the following assertion holds:  Let $\bu_0 \in \CD_{q,p}(\Omega)$,
and let  right members  
$\bff$, $g$, $\bg$, and $\bh$ for \eqref{lp}
satisfy the decay condition \eqref{cond:2}
and vanish for $t < 0$, then problem \eqref{lp}
with $T = \infty$ admits 
a unique solution 
$\bu \in L_p(\BR_+, W^2_q(\Omega)^N) \cap W^1_p(\BR_+, L_q(\Omega)^N)$
with some pressure term 
$\theta \in L_p(\BR_+, W^1_q(\Omega) + \hat W^1_{q,0}(\Omega))$ 
possessing the estimate:
\begin{equation}\label{est:2}
\|e^{\gamma t}\pd_t\bu\|_{L_p((0, T),L_q(\Omega))}
+ \|e^{\gamma t}\bu\|_{L_p((0, T), W^2_q(\Omega))}
\leq C\{ \bJ_{p,q}   
+ \delta(\Gamma)\sum_{\ell=1}^M
\Bigl(\int^T_0|e^{\gamma t}(\bu(\cdot, t), \bp_\ell)_\Omega|^p\,
dt\Bigr)^{1/p}\}
\end{equation}
for any $T > 0$ and $\gamma \leq \gamma_0$
with some constant $C$ independent of $T$.  
Here, 
$\delta(\Gamma)$ is a constant defined by  $\delta(\Gamma) = 1$
when $\Gamma=\emptyset$ and $\delta(\Gamma) = 0$ 
if $\Gamma \not= \emptyset$, and we have set 
$$\bJ_{p,q} =  \|\bu_0\|_{B^{2(1-1/p)}_{q,p}(\Omega)}
+ 
\|e^{\gamma t}(\bff, \tilde\Lambda_\gamma^{1/2}g, \pd_s\bg, 
\tilde\Lambda_\gamma^{1/2}\bh)
\|_{L_p(\BR, L_q(\Omega))} + 
\|e^{\gamma t}(g, \bh)\|_{L_p(\BR, W^1_q(\Omega))}.
$$
\end{thm}
\begin{proof}
We look for a solution $\bu$ of the form $\bu = \bv + \bw$, where  
$\bv$ and $\bw$ are  a solution to \eqref{gs*} and 
a solution to  problem:
\begin{gather}
\pd_t\bw - \DV\bT(\bw, \tilde\theta) = \lambda_1\bv, \quad
\dv \bw = 0 \quad  \text{in $\Omega\times(0, \infty)$}, 
\nonumber\\
\bT(\bw, \tilde\theta)\tilde\bn|_{S}  = 0, 
\quad \bw|_{\Gamma} = 0, 
\quad
\bw|_{t=0}  = 0 \quad\text{in $\Omega$}
\label{gs**}
\end{gather}
with some pressure term $\tilde \theta \in L_p(\BR_+, W^1_q(\Omega)
+ \hat W^1_{q,0}(\Omega))$, respectively.
 By Theorem \ref{thm:3.1}
\begin{equation}\label{est:3}
\|e^{\gamma t}\pd_t\bv\|_{L_p(\BR_+, L_q(\Omega))}
+ \|e^{\gamma t}\bv\|_{L_p(\BR_+, W^2_q(\Omega))}
\leq C\bJ_{p,q}.
\end{equation}
If $\Gamma\not=\emptyset$, setting 
$$\bw(\cdot, t) = \int^t_0 \BT(t-s)(\lambda_1\bv(\cdot, s))\,ds,$$
by Duhamel's principle we see that $\bw$ satisfies 
\eqref{gs**}.  Moreover, setting 
$$\bL_{q, \ba}(t) = \|\ba(\cdot, t)\|_{W^2_q(\Omega)}
+ \|\pd_t\ba(\cdot, t)\|_{L_q(\Omega)},$$
 by Theorem \ref{thm:exp-stab} 
and H\"older's inequality we have 
\begin{align*}
\bL_{q, \bw}(t) \leq C\int^t_0e^{-\gamma_0(t-s)}\bL_{q,\bv}(s)\,ds
\leq C(\gamma_0p')^{-1/{p'}}\Bigl(\int^t_0e^{-\gamma_0p(t-s)}
\bL_{q,\bv}(s)^p\,ds\Bigr)^{1/p}
\end{align*}
with some $\gamma_0 > 0$ for some positive constant $C$ independent of 
$t>0$, where $p' = p/(p-1)$.
 Thus, for $\gamma < \gamma_0$ we have 
\begin{align*}
\int^T_0(e^{\gamma t}\bL_{q,\bw}(t))^p\,dt
&\leq C(\gamma_0p')^{-p/{p'}}\int^T_0\bL_{q,\bv}(s)^p
\Bigl(\int^T_se^{-\gamma_0p(t-s)}e^{\gamma pt}\,dt\Bigr)\,ds\\
&=C(\gamma_0p')^{-p/{p'}}\int^T_0(e^{\gamma s}\bL_{q,\bv}(s))^p
\Bigl(\int^T_se^{-(\gamma_0-\gamma)p(t-s)}\,dt\Bigr)\,ds\\
&\leq C(\gamma_0p')^{-p/{p'}}((\gamma_0-\gamma)p)^{-1}\int^T_0(e^{\gamma s}
\bL_{q, \bv}(s))^p\,ds,
\end{align*}
which combined with \eqref{est:3} furnishes \eqref{est:2}
with $\delta(\Gamma) = 0$. 

Next, we consider the case of $\Gamma=\emptyset$.  Setting 
$\bz(x,t) = \lambda_1\bv(x,t) - \sum_{\ell=1}^M
(\lambda_1\bv(\cdot, t), \bp_\ell)_\Omega \bp_\ell(x)$, we have 
$(\bz(t), \bp_\ell)_\Omega = 0$ for $\ell=1, \ldots, M$ and $t > 0$.
Writing $\tilde\bw(t) = \int^t_0\BT(t-s)\bz(s)\,ds$, by 
Duhamel's principle, we see that $\tilde\bw$ satisfies \eqref{gs**} replacing 
$\lambda_1\bv$ by $\bz$.  Moreover, by Theorem \ref{thm:exp-stab}
and H\"older's inequality 
$$\bL_{q,\tilde\bw}(t) \leq C\int^t_0e^{-\gamma_0(t-s)}\bL_{q,\bz}(s)\,ds.$$
Thus, by \eqref{est:3} we have
\begin{equation}\label{est:4}
\int^T_0(e^{\gamma t}\bL_{q, \tilde\bw}(t))^p\,dt
\leq C_{\gamma_0, \gamma, p}\int^T_0(e^{\gamma t}\bL_{q, \bz}(t))^p\,dt
\leq C_{\gamma_0, \gamma, p}(\bJ_{p,q})^p.
\end{equation}
Setting $\bu = \bv + \bw$ with $\bw = \tilde\bw + \sum_{\ell=1}^M\int^t_0
(\lambda_1\bv(\cdot, s), \bp_\ell)_\Omega\bp_\ell$, we  see that 
$\bu$ satisfies \eqref{lp} with some pressure term
$\theta \in L_p(\BR_+, W^1_q(\Omega) + \hat W^1_{q,0}(\Omega))$,
because $\bD(\bp_\ell) = 0$ and $\dv \bp_\ell = 0$. 
Moreover, we have $\bD(\bu) = \bD(\bv) + \bD(\tilde\bw)$, so 
that by \eqref{est:3} and \eqref{est:4} we have
\begin{equation}\label{est:5}
\int^T_0\|\bD(\bu(\cdot, t))\|_{L_q(\Omega)}^p\,dt
\leq C(\bJ_{p,q})^p
\end{equation}
for any $T > 0$ with some constant $C$ independent of $T$. 
Since  
$$\|\ba\|_{W^1_q(\Omega)} \leq C(\|\bD(\ba)\|_{L_q(\Omega)}
+ \sum_{\ell=1}^M|(\ba, \bp_\ell)_\Omega|)
$$
for any $\ba \in W^1_q(\Omega)^N$ as follows from
 the usual contradiction argument
(cf. Duvaut and Lions \cite{DL}), by \eqref{est:5} 
\begin{equation}\label{est:6}
\int^T_0(e^{\gamma t}\|\bu(\cdot, t)\|_{W^1_q(\Omega)})^p\,dt
\leq C\{(\bJ_{p,q})^p 
+ \sum_{\ell-1}^M \int^T_0e^{\gamma pt}|(\bu(\cdot, t), 
\bp_\ell)_\Omega|^p\,dt\}.
\end{equation}
In addition, by \eqref{gre} with $\lambda=\lambda_0 + 1$, 
we have 
\begin{equation}\label{est:7}
\|\bu(t)\|_{W^2_q(\Omega)} \leq C\{\|\pd_t\bu(t)\|_{L_q(\Omega)}
+ \|\bu(t)\|_{L_q(\Omega)} 
+ \|(\bff(t), \bg(t))\|_{L_q(\Omega)}
+ \|(g(t), \bh(t))\|_{W^1_q(\Omega)}\}.
\end{equation}
Since $\pd_t\bu = \pd_t\bv + \pd_t\tilde\bw + \lambda_1\sum_{\ell=1}^M
(\bv(t), \bp_\ell)_\Omega \bp_\ell$, by \eqref{est:3}, 
\eqref{est:4}, \eqref{est:6} and \eqref{est:7} 
we have \eqref{est:2}, which completes the proof of Theorem \ref{thm:decay}.
\end{proof}
Finally, we prove the following theorem with help of Theorem
\ref{thm:decay}.
\begin{thm}\label{thm:decay*} Let $1 < p, q < \infty$, $N < r < \infty$ and 
$\max(q, q') \leq r$ $(q' = q/(q-1))$. Let $T$ be any positive number.
 Assume that 
$\Omega$ is a uniform $W^{2-1/r}_r$ domain and that $\Omega$ is bounded
in addition.  Then, there exists a positive constant $\gamma_0$ such that
for any $\bu_0 \in \CD_{p,q}(\Omega)$ and right members 
$\bff$, $g$, $\bg$ and $\bh$ with
\begin{alignat*}2
\bff &\in L_p((0, T), L_q(\Omega)^N), &\enskip
g &\in L_p((0, T), W^1_q(\Omega))\cap W^1_p((0, T), W^{-1}_q(\Omega)), \\
\bg &\in W^1_p((0, T), L_q(\Omega)^N)), &\enskip 
\bh &\in L_p((0, T), W^1_q(\Omega)^N) \cap W^1_p((0, T), W^{-1}_q(\Omega)^N)),
\end{alignat*}
satisfying the condition: $g|_{t=0} = 0$, $\bg|_{t=0} = 0$ and 
$\bh|_{t=0} = 0$, problem \eqref{lp} admits unique
solutions $\bu$ and $\theta$ with
$$\bu \in L_p((0, T), W^2_q(\Omega)^N) \cap W^1_p((0, T), L_q(\Omega)^N),
\enskip \theta \in L_p((0, T), W^1_q(\Omega) + \hat W^1_{q,0}(\Omega))$$
possessing the estimate:
\begin{align*}
&\|e^{\gamma t}\pd_t\bu\|_{L_p((0, t), L_q(\Omega))}
+ \|e^{\gamma t}\bu\|_{L_p((0, t), W^2_q(\Omega))}\\
&\quad \leq Ce^{2\gamma_0}\{\|\bu_0\|_{B^{2(1-1/p)}_{q,p}(\Omega)}
+ \delta(\Gamma)\sum_{\ell=1}^M\Bigl(\int^t_0 |
e^{\gamma s}(\bu(\cdot, s), \bp_\ell)_\Omega|^p\,ds
\Bigr)^{1/p} + \|e^{\gamma s}(\bff, \pd_t\bg)\|_{L_p((0, t), L_q(\Omega))}\\
&\quad + \|e^{\gamma s}(g, \bh)\|_{L_p((0, T), W^1_q(\Omega))}
+ \|e^{\gamma s}\pd_s[(1-\Delta)^{-1/2}(\iota g, \iota\bh)]
\|_{L_p((0, t), L_q(\BR^N))}\}
\end{align*}
for any $t \in (0, T]$ and 
$0 < \gamma \leq \gamma_0$ with some constant $C$ independent of 
$T$ and $\gamma$.  Here, $\delta(\Gamma)$ is the same number as in 
Theorem \ref{thm:decay}.
\end{thm}
\begin{proof}
Let $E_t$ be the same operator as in the proof of Theorem \ref{thm:linear*}.
Let $\phi(s)$ be a function in $C^\infty(\BR)$ such that $\phi(s) = 1$ 
for $s \leq 0$ and $\phi(s) = 0$ for $s \geq 1$ and set 
$\phi_t(s) = \phi(s-t)$.  Obviously, $\phi_t \in 
C^\infty(\BR)$, $\phi_t(s) = 1$ for $s \leq t$ and $\phi_t(s) = 0$
for $s \geq t+1$. Let $\bu^t=\bv$ and $\theta^t=\omega$ 
be solutions to the equations:
\begin{gather}
\pd_s\bv - \DV\bT(\bv, \omega) = \phi_tE_t\bff, \enskip
\dv\bv = \phi_tE_tg = \dv(\phi_tE_t\bg)
\quad\text{in $\Omega\times(0, \infty)$}, \nonumber\\
\bT(\bv, \omega)\tilde\bn|_S = \phi_tE_t\bh|_S, \enskip
\bv|_\Gamma = 0, \enskip \bv|_{s=0} = \bu_0 
\enskip\text{in $\Omega$}. \label{eq:3.1}
\end{gather}
Since $(\phi_tE_tf)(\cdot, s) = f(\cdot, s)$ for $s \in [0, T]$,
$\bu^t$ and $\theta^t$ solve problem \eqref{lp} for $s \in (0, t)$. 
And, by the uniqueness of solutions, $\bu^{t_1}(\cdot, s) =
\bu^{t_2}(\cdot, s)$ for $s \in [0, t_1]$ when $0 < t_1 < t_2 \leq T$.
By Theorem \ref{thm:decay}, 
\begin{align}
&\|e^{\gamma s}\bu^t\|_{L_p((0, t), W^2_q(\Omega))} 
+ \|e^{\gamma s}\pd_s\bu^t\|_{L_p((0, t), L_q(\Omega))} 
\leq C\{\|\bu_0\|_{B^{2(1-1/p)}_{q,p}(\Omega)} \nonumber\\
&\quad + \delta(\Gamma)\sum_{\ell=1}^M\Bigl(\int^t_0|e^{\gamma s}
|(\bu(\cdot, s), \bp_\ell)_\Omega|^p\,ds\Bigr)^{1/p} 
+ \|e^{\gamma s}(\phi_tE_t\bff, \pd_s(\phi_tE_t\bg))
\|_{L_p(\BR, L_q(\Omega))} \nonumber \\
&\qquad + \|e^{\gamma s}(\phi_tE_tg, \phi_tE_t\bh)
\|_{L_p(\BR, W^1_q(\Omega))}
 + \|e^{\gamma s}(\tilde\Lambda^{1/2}_\gamma(\phi_tE_tg), 
\tilde\Lambda^{1/2}_\gamma(\phi_TE_T\bh))
\|_{L_p(\BR, L_q(\Omega))}\}. \label{eq:3.2}
\end{align}
Let $X$ be $L_q(\Omega)$ or $W^1_q(\Omega)$.   
Using the change of variable: $2t-s = r$, we have
\begin{align*}
&\int^{2t}_te^{p\gamma s}\|(\phi_tE_tf)(\cdot, s)\|_X^p\,ds
= \int^{2t}_te^{p\gamma s}\phi_t(s)^p\|f(\cdot, 2t-s)\|_X^p\,ds \\
&\quad
\leq\int^t_{\max(0, t-1)} e^{2p\gamma(t-r)}e^{p\gamma r}
\|f(\cdot,r)\|_X^p\,dr
\leq e^{2p\gamma}\int^t_0e^{p\gamma r}\|f(\cdot, r)\|_X^p\,dr.
\end{align*}
Thus, noting that $\phi_tE_tf$ vanishes for $s \not\in [0, 2t]$,
we have 
\begin{equation}\label{eq:3.2*}
\|e^{\gamma s}\phi_tE_tf\|_{L_p(\BR, X)}
\leq e^{2\gamma}\|e^{\gamma s}f\|_{L_p((0, t), X)}.
\end{equation}
Noting \eqref{ext:2} and using \eqref{eq:3.2*}, we have
\begin{equation}\label{eq:3.3}\begin{split}
&\|e^{\gamma s}(\phi_tE_t\bff, \pd_s(\phi_tE_t\bg))
\|_{L_p(\BR, L_q(\Omega))}
\leq Ce^{2\gamma_0}\|e^{\gamma s}(\bff, \pd_t\bg)
\|_{L_p((0, t), L_q(\Omega))},  \\
&\|e^{\gamma s}(\phi_tE_tg, \pd_s(\phi_tE_t\bh))
\|_{L_p(\BR, W^1_q(\Omega))}
\leq Ce^{2\gamma_0}\|e^{\gamma s}(g, \bh)
\|_{L_p((0, t), W^1_q(\Omega))} 
\end{split}\end{equation}
for any $\gamma \in (0, \gamma_0]$ with some constant independent of 
$\gamma$, $t$ and $T$. 

In addition, applying the same argumentation as in the proof of 
the inequality \eqref{emb:2} in the appendix below, we have 
\begin{equation}\label{eq:3.3*}
\|e^{\gamma s}\tilde\Lambda_\gamma^{1/2}f\|_{L_p(\BR, L_q(\Omega))}
\leq C\{\|e^{\gamma s}\pd_s[(1-\Delta)^{-1/2}(\iota f)]
\|_{L_p(\BR, L_q(\BR^N))}
+ \|e^{\gamma s}f\|_{L_p(\BR, W^1_q(\Omega))}\},
\end{equation}
so that using \eqref{eq:3.2*} and \eqref{ext:2}, we have 
\begin{equation}\label{eq:3.4}\begin{split}
&\|e^{\gamma s}(\tilde\Lambda^{1/2}_\gamma(\phi_tE_tg), 
\tilde\Lambda^{1/2}_\gamma(\phi_tE_t\bh))
\|_{L_p(\BR, L_q(\Omega))} \\
&\quad
\leq Ce^{2\gamma_0}\{\|e^{\gamma s}\pd_s[(1-\Delta)^{-1/2}(\iota g, \iota\bh)]
\|_{L_p((0, t), L_q(\BR^N))}
+ \|e^{\gamma s}(g, \bh)\|_{L_p((0, t), W^1_q(\Omega))}\}.
\end{split}\end{equation}
Setting $\bu = \bu^T$ and $\theta = \theta^T$ 
and combining \eqref{eq:3.2}, \eqref{eq:3.3} and \eqref{eq:3.4},
we have Theorem \ref{thm:decay*}. 
\end{proof}
\section{A proof of a global in time unique existence theorem}

In this section, we prove Theorem \ref{main:global}, so that we assume that
$\Omega$ is bounded in addition. Let $T_0$ be a positive number such that 
for any initial data $\bv_0 \in \CD_{q,p}(\Omega)$ with 
$\|\bv_0\|_{B^{2(1-1/p)}_{q,p}(\Omega)} \leq 1$, problem \eqref{p1}
admits a unique solution $\bu \in \CS^{1,2}_{p,q}(0, T_0)$ satisfying
\eqref{small:1}.   Here and in the following, we set  
\begin{align*}
\CS^{1,2}_{p,q}(a,b) &= W^1_p((a, b), L_q(\Omega)^N) 
\cap L_p((a, b), W^2_q(\Omega)^N), \\
\BI_\bv(a, b) &= \|e^{\gamma t}\pd_t\bv\|_{L_p((a, b), L_q(\Omega))}
+ \|e^{\gamma t}\bv\|_{L_p((a, b), W^2_q(\Omega))}
\end{align*}
for any $a$, $b$ satisfying  $0 \leq a < b \leq \infty$
for the notational simplicity, where 
$\gamma$ is a fixed positive number for which 
Theorem \ref{thm:exp-stab} and 
Theorem \ref{thm:decay} hold.
By Theorem \ref{main:loc}, 
such $T_0 > 0$ exists. 

Let $\epsilon$ be a small positive number $\leq 1$
that is determined later 
and we assume that $\bv_0 \in \CD_{q,p}(\Omega)$ and 
$\|\bv_0\|_{B^{2(1-1/p)}_{q,p}(\Omega)} \leq \epsilon$. Let $T$
be a positive number such that  
problem \eqref{p1} admits a solution $\bu \in \CS^{1,2}_{p,q}
(0, T)$ that satisfies \eqref{small:1}. 
Since $\|\bv_0\|_{B^{2(1-1/p)}_{q,p}(\Omega)} \leq \epsilon \leq 1$, 
we have $T \geq T_0$.  The main step 
is to prove that there exist constants $\epsilon_0>0$
and $M_6$ independet of $\epsilon$ and $T$ such that
\begin{equation}\label{global:est1}
\BI_\bv(0, t) \leq M_6(\epsilon + \BI_\bv(0, t)^2)
\end{equation}
for any $t \in (0, T]$ provided that $0 < \epsilon \leq \epsilon_0$. 

In fact, let $r_\pm(\epsilon)$ be two roots of the quadratic equation:
$M_6(\epsilon + x^2) - x=0$, that is $r_\pm(\epsilon) = 
(2M_6)^{-1} \pm \sqrt{(2M_6)^{-2}-\epsilon}$. We find a small posivite
number $\epsilon_1 > 0$ such that $0 < r_-(\epsilon) < r_+(\epsilon)$
whenever $0 < \epsilon \leq \epsilon_1$.  In this case, 
$r_-(\epsilon ) = M_6\epsilon + O(\epsilon^2)$ as $\epsilon\to0+0$.
Since $\BI_\bu(0, t) \to 0$ as $t\to0$ and $\BI_\bu(0,t)$ is a
continuous function with respect to $t$, by 
\eqref{global:est1} we have $\BI_\bu(0, t) 
\leq r_-(\epsilon)$ for any $t \in (0, T]$, especially $\BI_\bu(0, T)
\leq r_-(\epsilon)$. To prove 
\begin{equation}\label{global:est2}
\sup_{t \in (0, T)} \|\bu(\cdot, t)\|_{B^{2(1-1/p)}_{q,p}(\Omega)}
\leq M_7(\bI_\bu(0, T) + e^{-\gamma T}
\|\bv_0\|_{B^{2(1-1/p)}_{q,p}(\Omega)}),
\end{equation}
we take $\bZ \in \CS^{1,2}_{p,q}(0, \infty)$ which solves \eqref{5.3}
with some pressure term $\Pi$.  By Corollary \ref{exp-stab2}
we have
$$\|e^{\gamma t}\pd_t\bZ\|_{L_p(\BR_+, L_q(\Omega))}
+ \|e^{\gamma t}\bZ\|_{L_p(\BR_+, W^2_q(\Omega))} \leq 
M_8\|\bv_0\|_{B^{2(1-1/p)}_{q,p}(\Omega)}$$
with some constant $M_8$ independent of $\epsilon$ and $T$,
because $\bv_0$ satisfies \eqref{orth:1} when $\Gamma = \emptyset$.  
Employing the same argument as in the proof of \eqref{eq:2.4}, 
we have \eqref{global:est2}. 

Since we may assume that $r_-(\epsilon) \leq 2M_0$,
by  \eqref{global:est2} we have
$\|\bu(\cdot, T-0)\|_{B^{2(1-1/p)}_{q,p}(\Omega)} \leq M_7(2M_0+1)\epsilon$,
because $e^{-\gamma T} \leq 1$.  Choose $\epsilon$ so small that 
$M_7(2M_0 + 1)\epsilon \leq 1$.  By Theorem \ref{main:loc},
there exists a unique solution $\bu' \in \CS^{1,2}_{p,q}(T,  T+T_0)$
of the equations:
\begin{gather*}
\pd_t\bu' - \DV\bT(\bu', \theta') = \bF(\bu'),
\enskip \dv\bu' = G(\bu') = \dv\bG(\bu') \quad\text{in 
$\Omega\times(T, T+T_0)$}, \\
\bT(\bu', \theta')\bn|_S = \bH(\bu')\bn|_S,
\enskip \bu'|_\Gamma=0, \enskip \bu'|_{t=T+0} = \bu|_{t=T-0},
\end{gather*}
with some pressure term $\theta'$. Choosing $T_0$ smaller if necessary,
we may assume that 
$$\int^{T+T_0}_T\|\nabla\bu'(\cdot, t)
\|_{L_\infty(\Omega)}\,dt \leq \sigma/2.$$
Since 
$\int^T_0\|\nabla\bu(\cdot, t)
\|_{L_\infty(\Omega)}\,dt \leq M_9\BI_\bu(0, T)$ with some 
constant $M_9$ independent of $\epsilon$ and $T$ as follows from 
\eqref{5.1}, we choose $\epsilon$ so small that $M_9\epsilon <\sigma/2$,
so that $\int^T_0\|\nabla\bu(\cdot, t)
\|_{L_\infty(\Omega)}\,dt \leq \sigma/2$.  If we define
$\bu''$ by $\bu''(\cdot, t) = \bu(\cdot, t)$ for $0 \leq t \leq T$
and $\bu''(\cdot, t) = \bu'(\cdot, t)$ for $T \leq t \leq T+T_0$,
then $\bu''$ satisfies the equations \eqref{p1}
for $t \in (0, T+T_0)$ with some pressure term 
$\theta'$ and the condition: $\int^{T+T_0}_0
\|\bu''(\cdot, t)\|_{L_\infty(\Omega)}\,dt \leq \sigma$.
Thus, by \eqref{global:est1} $\BI_{\bu''}(0, T+T_0)
\leq M_6(\epsilon + \BI_{\bu''}(0, T+T_0)^2)$.  Repeating this 
argument, we can prolong $\bu$ to any time interval $(0, T)$
with $\BI_{\bu}(0, T) \leq r_-(\epsilon)$, which completes 
the existence of solution $\bu$ globally defined in time with
$\BI_\bu(0, \infty) \leq r_-(\epsilon)$.  The uniquness follows
from the same argumentations as in the proof of Theorem \ref{main:loc}
with small $\epsilon > 0$ instead of small $T>0$. Therefore,
our task is to prove \eqref{global:est1}.

Applying Theorem \ref{thm:decay} to problem \eqref{p1}, 
we have 
\begin{equation}\label{est:4.1}
\BI_\bu(0, t)
\leq 
C\{\|\bv_0\|_{B^{2(1-1/p)}_{q,p}(\Omega)}
+ \delta(\Gamma)\sum_{\ell=1}^M
\Bigl(\int^t_0|e^{\gamma s}(\bu(\cdot, s), \bp_\ell)_\Omega||^p\,ds
\Bigr)^{1/p} + \BK_\bu(0, t)\}
\end{equation}
for any $t \in (0, T]$ with 
\begin{align*}\BK_\bu(0, t) &= 
 \|e^{\gamma s}(\bF(\bu), \pd_s\bG(\bu))\|_{L_p((0, t), L_q(\Omega))}
+ \|e^{\gamma s}(G(\bu), \bH(\bu)\tilde\bn)\|_{L_p((0, t), W^1_q(\Omega))}\\
&+ \|e^{\gamma s}\pd_s[(1-\Delta)^{-1/2}(\iota G(\bu), 
\iota \bH(\bu)\tilde\bn))\|_{L_p((0, t), L_q(\BR^N))}.
\end{align*}
Here and in the following, $C$ denotes a generic constant independent of 
$\epsilon$, $t \in (0, T]$ and $T$. 

When $\Gamma=\emptyset$, $\delta(\Gamma) = 1$, so that we have to estimate
$\int^t_0|e^{\gamma t}(\bu(\cdot, s), \bp_\ell)_\Omega|)^p\,ds$.
Recalling Remark \ref{rem:5} \thetag1, $\bv(x, t) 
= \bu(\bX_\bu^{-1}(x, t), t)$ satisfies the equation
\eqref{1.1} with \eqref{1.2}, where $\bX^{-1}(x, t)$ denotes the 
inverse map of the correspondence: $x = \xi + \int^t_0\bu(\xi, s)\,ds
= \bX_\bu(\xi, t)$.  Since 
$\frac{d}{dt}\int_{\Omega_t} \bv(x, t)\bp_\ell(x)\,dx = 0$, by \eqref{orth:1}
we have $\int_{\Omega_t}\bv(x, t)\bp_\ell(x)\,dx = 0$, so that
$$\int_\Omega \bu(\xi, s))\bp_\ell(\xi+ \int^s_0 \bu(\xi, r)\,dr)\,d\xi=0,
$$
which combined with \eqref{5.1} and H\"older's inequality furnishes that 
\begin{align*}
|(\bu(\cdot, s), \bp_\ell)_\Omega|& \leq C\|\bu(\cdot, s)\|_{L_q(\Omega)}
\int^s_0\|\bu(\cdot, r)\|_{W^1_q(\Omega)}\,dr \\
& \leq C\|\bu(\cdot, s)\|_{L_q(\Omega)}
\Bigl(\int^s_0e^{-p'\gamma r}\,ds\Bigr)^{1/{p'}}
\Bigl(\int^t_0 (e^{\gamma r}\|\bu(\cdot, r)\|_{W^1_q(\Omega)})^p\,dr
\Bigr)^{1/p}.
\end{align*}
Thus, we have 
\begin{equation}\label{est:4.3}
\delta(\Gamma) \sum_{\ell=1}^M \Bigl(\int^t_0 |e^{\gamma s}
(\bu(\cdot, s), \bp_\ell)_\Omega|^p\,ds\Bigr)^{1/p}
\leq C\BI_\bu(0, t)^2.
\end{equation}

From now on, we estimate $\BK_\bu(0,t)$. By H\"older's inequality
we have 
\begin{equation}\label{eq:4.2}
\sup_{s\in(0, t)}\|\int^s_0\nabla\bu(\cdot,r)\,ds\|_{W^1_q(\Omega)}
\leq \Bigl(\int^s_0(e^{\gamma r}\|\bu(\cdot, r)\|_{W^2_q(\Omega)}
)^p\,dr\Bigr)^{1/p}\Bigl(\int^s_0e^{-\gamma rp'}\,dr\Bigr)^{1/{p'}}
\leq C\BI_\bu(0,t).
\end{equation}
Since \eqref{small:1} holds and since we may assume that
$\sigma\leq 1$, by \eqref{null:1}, \eqref{eq:4.2} and \eqref{5.1}
\begin{align}
&\sup_{s \in (0, t)}\|\int^s_0\nabla\bu(\cdot, r)\,dr\|_{L_\infty(\Omega)}
\leq 1, \enskip 
\sup_{s\in(0, t)} \|\bV_i(\int^s_0 \nabla\bu(\cdot, r)\,dr
\|_{W^1_q(\Omega)} \leq C\,\BI_\bu(0, t), \nonumber \\
&\sup_{s\in(0, t)} \|\nabla \bW(\int^s_0 \nabla\bu(\cdot, r)\,dr
\|_{L_q(\Omega)} \leq C\,\BI_\bu(0, t), \label{eq:4.3}
\end{align}
where $i = 1,2,4,5$ and $6$ and $\bW = \bW(\bK)$ is any matrix of 
polynomials with respect to $\bK$.  
By \eqref{5.1}, \eqref{ext-normal}, \eqref{nonlinear:1} and 
\eqref{eq:4.3}, we have 
\begin{equation}\label{eq:4.4}
\|\bF(\bu)\|_{L_p((0, t), L_q(\omega))} \leq C\BI_\bu(0, t)^2, \quad
\|(G(\bu), \bH(\bu)\tilde\bn)\|_{L_p((0, t), W^1_q(\Omega))}
\leq C\BI_\bu(0, t)^2.
\end{equation}

Since $B^{2(1-1/p)}_{q,p}(\Omega) \subset W^1_q(\Omega)$ as follows
from the assumption: $2 < p , \infty$, by \eqref{global:est2}
we have
\begin{align}
&\sup_{s\in(0,t)} \|\bu(s, \cdot)\|_{W^1_q(\Omega)}
\leq C(\BI_\bu(0, t) + \epsilon), \enskip 
\sup_{s\in(0,t)} \|\pd_s\bW(\int^s_0\nabla \bu(\cdot, r)\,
dr)\|_{L_q(\Omega)}
\leq C(\BI_\bu(0, t) + \epsilon).
\label{eq:4.4*}
\end{align}
Thus, by \eqref{5.1}, \eqref{eq:4.3} and \eqref{eq:4.4*} 
\begin{equation}\label{eq:4.5}
\|\pd_s\bG(\bu)\|_{L_p((0, t), L_q(\Omega))}
\leq C(\bI_\bu(0, t)^2 + (\BI_\bu(0,T)+\epsilon)\BI_\bu(0,T))
\leq 2C(\BI_\bu(\bu)^2 + \epsilon),
\end{equation}
because $0 < \epsilon \leq 1$.

To estimate $\|e^{\gamma t}(\tilde\Lambda^{1/2}_\gamma g_\bu, 
\tilde\Lambda^{1/2}_\gamma (\bh_\bu\tilde\bn))\|_{L_p(\BR, L_q(\Omega))}$,
we use the following lemma which can be proved 
in the same manner as in the
proof of Lemma \ref{lem:est1}. 
\begin{lem}\label{lem:4.1} Let $1 < p, q < \infty$, $N < r, q < \infty$
and let 
$\Omega$ be a uniform $W^{2-1/r}_r$. Let $\iota$ be the 
extension map satisfying the properties \thetag{e-1} and 
\thetag{e-2}.  Then 
\begin{multline*}
\|e^{\gamma s}\tilde\Lambda_\gamma^{1/2}((\nabla f)g)
\|_{L_p((0,t), L_q(\Omega))}^p \\
\leq C\Bigl\{\int^t_0(e^{-\gamma s}\|\pd_sf(\cdot, s)\|_{L_q(\Omega)}
\|g(\cdot, s)\|_{W^1_q(\Omega)})^p\,ds 
+ \int^t_0(e^{\gamma s}
\|\nabla f(\cdot, s)\|_{L_q(\Omega)}
\|\pd_tg(\cdot, s)\|_{L_q(\Omega)})^p\,ds\Bigr\}.
\end{multline*}
\end{lem}
Applying Lemma \ref{lem:4.1} and using  \eqref{ext-normal},
\eqref{eq:4.3} and \eqref{eq:4.4}, we have 
\begin{equation}\label{est:4.5}\begin{split}
&\|e^{\gamma s}\pd_s[(1-\Delta)^{-1/2}(\iota G(\bu), 
\iota \bH(\bu)\tilde\bn)]\|_{L_p((0, t), L_q(\BR^N))} \\
&\quad
\leq C(\BI_\bu(0,t)^2 + (\BI_\bu(0,T) + \epsilon)\BI_\bu(0,T)) \leq 
2(\BI_\bu(0, T)^2 + \epsilon),
\end{split}\end{equation}
which combined with \eqref{est:4.1}, \eqref{est:4.3}, \eqref{eq:4.4}
and \eqref{eq:4.5} furnishes \eqref{global:est1}.
This completes the proof of Theorem \ref{main:global}.

\appendix
\section{A proof of the inequality \eqref{emb:2}}
First, we prove the inequality \eqref{emb:2} in case of $\Omega = \BR^N$.
\begin{lem}\label{ap:1} Let $1 < p, q < \infty$ and set 
$(1-\Delta)^sf - \CF^{-1}_\xi[(1+|\xi|^2)^{s/2}\hat f(\xi)]$
for $f \in \CS'(\BR^N)$ and $s \in \BR$.  Here, $\hat f$ denotes
the Fourier transform of $f$, $\CF^{-1}_\xi$ the inverse Fourier
transform and $\CS'(\BR^N)$ the space of tempered distributions
on $\BR^N$ in the sense of L. Schwartz. Then, we have
\begin{align*}
&\|e^{-\gamma t}\Lambda_\gamma^{1/2}f\|_{L_p(\BR, L_q(\BR^N))}\\
&\quad\leq C\{\|e^{-\gamma t}f\|_{L_p(\BR, L_q(\BR^N))}
+ \|e^{-\gamma t}(1 - \Delta)^{-1/2}\pd_tf\|_{L_p(\BR, L_q(\BR^N))}
+ \|e^{-\gamma t}(1-\Delta)^{1/2}f\|_{L_p(\BR, L_q(\BR^N))}\}.
\end{align*}
\end{lem}
\begin{proof}
The idea of our proof here is the same as in the proof of Proposition 
2.8 in \cite{SS2}. Let $\varphi_0(t)$ be a function in 
$C^\infty(\BR)$ such that $\varphi_0(t) = 1$ for $|t| \leq 1$ and 
$\varphi_0(t) = 0$ for $|t| \geq 2$, and set $\varphi_\infty(t) 
= 1 - \varphi_0(t)$.  We define functions $A_j(\xi, \lambda)$ ($j = 1, 2$)
by 
\begin{align*}
A_1(\xi, \lambda) &= \varphi_\infty(\tau)\varphi_0
\Bigl(\frac{(1+|\xi|^2)^{1/2}}{|\lambda|}\Bigr)
\frac{(1+|\xi|^2)^{1/4}\lambda^{1/2}}{\lambda}, 
\\
A_2(\xi, \lambda) &= \varphi_\infty(\tau)
\varphi_\infty
\Bigl(\frac{(1+|\xi|^2)^{1/2}}{|\lambda|}\Bigr)
\frac{\lambda^{1/2}}{(1 + |\xi|^2)^{1/4}}.
\end{align*}
We have
$$|\pd_\tau^\ell\pd_\xi^\alpha A_j(\xi, \lambda)|
\leq C_{\ell, \alpha}|\tau|^{-\ell}|\xi|^{-|\alpha|}
\quad(\lambda = i\tau + \gamma,\enskip j=1, 2)
$$
for any $\ell \in \BN_0$ and $\alpha \in \BN_0^N$, 
$\xi \in \BR^N\setminus\{0\}$, $\tau, \gamma \in \BR\setminus\{0\}$
with some constant $C_{\ell, \alpha}$ depending solely on 
$\ell$ and $\alpha$. Set $A_j(\lambda, D_x)f = 
\CF^{-1}_\xi[A_j(\xi, \lambda)\hat f(\xi)]$ for any $f \in \CS'(\BR^N)$,
and then by Theorem 3.3 in \cite{ES} we know that the sets
$\{\tau^k\pd_\tau^kA_j(\lambda, D_x) \mid \tau \in \BR\setminus\{0\}\}$
are $\CR$-bounded families in $\CL(L_q(\BR^N))$ and their
$\CR$-bounded are less than $C_{q,N}\max_{|\alpha| \leq N+2} 
C_{k, \alpha}$ for $k = 0, 1$ and $j = 1, 2$, where 
$\CL(L_q(\BR^N))$ is the set of all bounded linear operators on 
$L_q(\BR^N)$.  Therefore, by Weis's operator valued Fourier
multiplier theorem \cite{Weis} we have
\begin{equation}\label{Ap:1}
\|e^{-\gamma t}A_j(\pd_t, D_x)F\|_{L_p(\BR, L_q(\BR^N))}
\leq C\|e^{-\gamma t}F\|_{L_p(\BR, L_q(\Omega))}
\quad(\gamma \not=0).
\end{equation}
Here, the operators $A_j(\pd_t, D_x)$ are defined by 
$$A_j(\pd_t, D_x)F = \CL^{-1}_\lambda[A_j(\lambda, D_x)
\CL[F](\lambda, \cdot)](t).$$
Dividing $\lambda^{1/2}$ into the following three parts:
$$\lambda^{1/2} = \varphi_0(\tau)\lambda^{1/2}
+ \varphi_\infty(\tau)A_1(\xi, \lambda)\frac{\lambda}{(1 + |\xi|^2)^{1/4}}
+ \varphi_\infty(\tau)A_2(\xi, \lambda)(1 + |\xi|^2)^{1/4},$$
and using \eqref{Ap:1} and Bourgain's Fourier multiplier theorem
\cite{Bourgain}, we have Lemma \ref{Ap:1}.
\end{proof}
{\bf Proof of the inequality \eqref{emb:2}.}~ To prove the lemma, we use 
the exitension map $E$ having the following properties:
\begin{itemize}
\item[\thetag{Ex-1}]~ $\|Ef\|_{L_q(\BR^N)} 
\leq C_{q, \Omega}\|f\|_{L_q(\Omega)}$
\item[\thetag{Ex-2}]~ $\|(1 - \Delta)^{1/2}Ef\|_{L_q(\BR^N)}
\leq C_{q, \Omega}\|f\|_{W^1_q(\Omega)}$
\item[\thetag{Ex-3}]~ $\|(1 - \Delta)^{-1/2}E(\nabla f)\|_{L_q(\BR^N)}
\leq C_{q, \Omega}\|f\|_{L_q(\Omega)}$.
\end{itemize}
Such extension map can be constructed under the assumption that
$N < q, r < \infty$. By Lemma \ref{Ap:1}, we have 
\begin{align*}
&\|e^{-\gamma t}\Lambda_\gamma^{1/2}((\nabla f)g)\|_{L_p(\BR, L_q(\Omega))}
\leq \|e^{-\gamma t}\Lambda_\gamma^{1/2}E((\nabla f)g)
\|_{L_p(\BR, L_q(\BR^N))}\\
&\quad\leq \|e^{-\gamma t}E((\nabla f)g)\|_{L_p(\BR, L_q(\Omega))} 
+ \|e^{-\gamma t}(1-\Delta)^{-1/2}\pd_tE((\nabla f)g)\|_{L_p(\BR, L_q(\Omega))}
\\
&\quad\quad
 + \|e^{-\gamma t}(1-\Delta)^{1/2}((\nabla f)g)\|_{L_p(\BR, L_q(\Omega))}.
\end{align*}
Using the identity: $\pd_t((\nabla f)g) = \nabla(\pd_tf\cdot g)
- \pd_tf(\nabla g) + (\nabla f)\pd_tg$ and \thetag{Ex-3}, we have
\begin{align*}
&\|e^{-\gamma t}(1-\Delta)^{-1/2}\pd_tE((\nabla f)g)\|_{L_p(\BR, L_q(\Omega))}
\leq C\|e^{-\gamma t}(\pd_tf\cdot g)\|_{L_p(\BR, L_q(\Omega))} \\
&\quad 
+ \|e^{-\gamma t}(1-\Delta)^{-1/2}E(\pd_tf(\nabla g))\|_{L_p(\BR, L_q(\BR^N))}
+ \|e^{-\gamma t}(1-\Delta)^{-1/2}E((\nabla f)\pd_tg)\|_{L_p(\BR, L_q(\BR^N)}.
\end{align*}
By \eqref{5.1} and \thetag{Ex-1}, we have
\begin{align*}
\|e^{-\gamma t}(\pd_tf\cdot g)\|_{L_p(\BR, L_q(\Omega))}^p
&\leq C\int^\infty_{-\infty}(e^{-\gamma t}
\|\pd_tf(\cdot, t)\|_{L_q(\Omega)}\|g(\cdot, t)\|_{W^1_q(\Omega)})^p\,dt, \\
\|e^{-\gamma t}E((\nabla f)g)\|_{L_p(\BR, L_q(\Omega))}^p
&\leq C\int^\infty_{-\infty}(e^{-\gamma t}
\|\nabla f(\cdot, t)\|_{L_q(\Omega)}\|g(\cdot, t)\|_{W^1_q(\Omega)})^p\,dt.
\end{align*}
To estimate other terms, we use the inequality:
\begin{equation}\label{ap.2}
\|(1-\Delta)^{-1/2}E(fg)\|_{L_q(\BR^N)}
\leq C\|f\|_{L_q(\Omega)}\|g\|_{L_q(\Omega)}.
\end{equation}
In fact, for any $\varphi \in C^\infty_0(\BR^N)$, we observe that
\begin{align*}
|((1-\Delta)^{-1/2}E(f, g), \varphi)_{\BR^N}|
&\leq \|E(fg)\|_{L_{q/2}(\BR^N)}\|(1-\Delta)^{-1/2}\varphi\|_{L_s(\BR^N)}\\
&\leq C\|f\|_{L_q(\Omega)}\|g\|_{L_q(\Omega)}
\|(1-\Delta)^{-1/2}\varphi\|_{L_s(\BR^N)},
\end{align*}
where $s$ is an index such that $1/s + 2/q = 1$.  Since $2 \leq N < q < \infty$,we can choose such $s$ with $1 < s < \infty$.  Since 
$N(1/q' - 1/s) = N(1-1/q-1/s) = N/q <1$, we have 
$\|(1-\Delta)^{-1/2}\varphi\|_{L_s(\BR^N)} \leq C\|(1-\Delta)^{-1/2}\varphi
\|_{W^1_{q'}(\BR^N)}$, so that we have \eqref{ap.2}.

By \eqref{ap.2} we have 
\begin{align*}
\|e^{-\gamma t}(1 - \Delta)^{-1/2}E(\pd_tf(\nabla g))\|_{L_p(\BR, L_q(\BR^N))}^p
\leq C\int^\infty_{-\infty}(e^{-\gamma t}
\|\pd_tf(\cdot, t)\|_{L_q(\Omega)}\|\nabla g(\cdot, t)\|_{L_q(\Omega)})^p\,dt,\\
\|e^{-\gamma t}(1 - \Delta)^{-1/2}E((\nabla f)\pd_tg))\|_{L_p(\BR, L_q(\BR^N))}^p
\leq C\int^\infty_{-\infty}(e^{-\gamma t}
\|\nabla f(\cdot, t)\|_{L_q(\Omega)}
\|\pd_tg(\cdot, t)\|_{L_q(\Omega)})^p\,dt.
\end{align*}
By \eqref{5.1} and \thetag{Ex-2} we have
\begin{align*}
\|e^{-\gamma t}(1 - \Delta)^{1/2}E((\nabla f)g)\|_{L_p(\BR, L_q(\BR^N))}^p
&\leq
C\|e^{-\gamma t}((\nabla f)g)\|_{L_p(\BR, W^1_q(\BR^N))}^p \\
&\leq C\int^\infty_{-\infty}(e^{-\gamma t}
\|\nabla f(\cdot, t)\|_{W^1_q(\Omega)}
\|\nabla g(\cdot, t)\|_{W^1_q(\Omega)})^p\,dt.
\end{align*}
This completes the proof of the inequality \eqref{emb:2}. 

\end{document}